\documentclass[a4paper,11pt]{amsart}
\usepackage[margin=1.2in]{geometry}

\usepackage{mathpazo,amsmath,amssymb,amscd,color,xcolor,mathrsfs}
\usepackage{tikz}
\usepackage{float}
\usepackage{comment}
\definecolor{green}{RGB}{0,144,0}
\definecolor{bluegreen}{RGB}{17,100,180}
\usepackage[linkcolor = bluegreen, citecolor = bluegreen, colorlinks = True]{hyperref}
\numberwithin{equation}{section}
\numberwithin{figure}{section}
\numberwithin{table}{section}
\usepackage[all]{hypcap}
\usepackage{mathtools}

\newtheorem{theorem}{Theorem}[section]
\newtheorem*{theorem*}{Theorem}

\newtheorem{lemma}[theorem]{Lemma}
\newtheorem{corollary}[theorem]{Corollary}

\newtheorem{question}[theorem]{Question}

\newtheorem{definition}[theorem]{Definition}

\DeclareMathOperator{\M}{\text{M}}

\title[Pseudo-Anosovs from mapping tori]{Pseudo-Anosovs from the perspective of their mapping tori}

\author{Tarik Aougab}
\address{Department of Mathematics, Haverford College, Haverford, PA 19041, USA}
\curraddr{}
\email{taougab@haverford.edu}

\date{}

\begin{document}

\begin{abstract} In this chapter, we outline some of the many combinatorial tools developed over the past three decades for studying a pseudo-Anosov diffeomorphism of a surface by analyzing the geometry of its mapping torus. We begin with an overview of the various simplicial complexes associated with a surface (such as the curve, arc, and pants complexes) and explain how to relate the dynamics of the action of a given pseudo-Anosov on any one of these complexes to the dynamics of the diffeomorphism itself, or to the hyperbolic geometry of its mapping torus. We next cover some of the more modern features of the theory by discussing various analogs of pseudo-Anosov diffeomorphisms on surfaces of infinite type. We conclude with a description of original work-- due jointly to the author with Dave Futer and Sam Taylor-- that relates the action of a pseudo-Anosov on the curve complex to the minimum number of fixed points for any map in the corresponding isotopy class. The paper is written in as accessible a way as possible while assuming only the bare minimum in background. The hope is to informally convey to the reader some of the main ideas and strategies in the area.

\end{abstract}

\maketitle


\textbf{Relevant AMS codes:} 57K20, 57K32, 57M07, 57M60


\section{Introduction}\label{sec:intro}

Thurston's beautiful and timeless work on the theory of surface diffeomorphisms has led to an explosion of new mathematics, tied to closely related areas in group theory, dynamics, and geometry \cite{Thurston1}. The famous Nielsen--Thurston trichotomy partitions isotopy classes of diffeomorphisms into three types, and the one of those with by far the richest and most interesting dynamical properties are the \textit{pseudo-Anosovs}\index{pseudo-Anosov}: those which have a representative that leave invariant a pair of transverse singular foliations and which stretch along one and contract along the other. Over the last thirty years, tremendous progress has been made in the quest to fully understand the various invariants associated with a pseudo-Anosov class (\cite{MasurMinskyII}, \cite{Brock2}, \cite{MT}, \cite{FS}, \cite{TL}). The humble goal of this expository piece is to detail some of this progress, particularly the work that acts as a bridge between tools with a more combinatorial flavor, and those that are more geometric.

Central to this story is the famous \textit{complex of curves}\index{curve complex}-- as well as its relatives, the \textit{arc complex}\index{arc complex} and the \textit{pants complex}\index{pants complex}-- associated with an orientable surface $S$ of negative Euler characteristic. The mapping class group\index{mapping class group} acts by simplicial automorphisms on the curve complex and its relatives. Groundbreaking work of Brock, Minsky, and Masur--Minsky (\cite{Brock1}, \cite{Brock2}, \cite{Minsky}, \cite{MasurMinskyI}, \cite{MasurMinskyII}) which eventually lead to Brock--Canary--Minsky's celebrated resolution of Thurston's \textit{ending lamination conjecture} \index{ending laminations} \cite{BCM}, allows one to read off information about a pseudo-Anosov from its action on these complexes.

For example, letting $\phi$ denote a pseudo-Anosov mapping class,  Thurston's hyperbolization theorems imply that the mapping torus $M_{\phi}$ admits a complete hyperbolic metric \cite{Thurston3}. In turn, Mostow--Prasad rigidity\index{Mostow-Prasad rigidity} implies that geometric properties of this hyperbolic metric-- for example its volume, or the length of the shortest closed geodesic-- can be taken to be invariants of the pseudo-Anosov $\phi$ itself: if a pair of pseudo-Anosov diffeomorphisms are isotopic, they must agree on all such geometric quantities. 

The purpose of this survey is to highlight some key examples of this bridge between combinatorics, dynamics, and geometry, and then to detail some of the more modern and recent contributions to this story. After going through some basic definitions and preliminaries in Section \ref{prelim}, we will begin in earnest in Section \ref{geometry from combinatorics} with some detailed descriptions of results that allow one to read off information about the geometry of a hyperbolic mapping torus from the action of its monodromy on a combinatorial complex\index{simplicial complex} associated with its fiber surface.

Agol's proof of the Virtual Fibering and Virtual Haken conjectures in some sense completes a long-standing program to relate the structure of hyperbolic $3$-manifolds to the dynamics of surface diffeomorphisms \cite{Agol3}. But what remains completely mysterious is how to generalize any of this to the world of \textit{infinite type surfaces}, those whose fundamental groups are not finitely generated. Much work has been done in this area in the past half--decade-- we refer the reader to Aramayona-Vlamis' excellent expository piece for more details \cite{AV}. Perhaps unsurprisingly, many of the recent results on the mapping class group of an infinite type surface are \textit{negative}, in the sense that they show that a particular property of a finite type mapping class group no longer holds in this much wilder setting. What is far more interesting and useful are those rarer \textit{positive} results, ones that show how techniques in the finite type world still apply in the infinite setting. In Section \ref{infinite type}, we will document some of this recent work, in particular highlighting results of Field-Kent-Kim-Leininger-Loving  (\cite{FKeLL}, \cite{FKLL}) and Whitfield \cite{Whitfield}.

We then move to an exposition of new work of the author-- joint with Futer and Taylor \cite{AFT} -- regarding the fixed points\index{fixed points!pseudo-Anosov} of a pseudo-Anosov (more precisely, the minimum number of fixed points taken over all representatives in a pseudo-Anosov isotopy class). We show that, under mild --and necessary-- assumptions on the pseudo-Anosov, the logarithm of the number of fixed points is proportional to its Teichm{\"u}ller dilatation. Without these minor assumptions, we obtain another result that allows one to estimate the number of fixed points from various combinatorial data, and we explain how this more general result is optimal. 

Finally, we conclude the piece with some open questions.

\subsection{How to read this piece.}
Through the use of a detailed preliminaries section (Section \ref{prelim}), we do our best to make the piece as self-contained as possible. We will assume only a general background on the topology of $2$- and $3$-dimensional manifolds, and some basic properties of hyperbolic geometry. Therefore, by reading Section \ref{prelim}, a mathematically literate reader who is unfamiliar with the main tools and techniques in this area should be able to pick up the minimum background required for appreciating the results that follow. 

Each subsection in Section \ref{prelim} begins with a ``top-level summary'', which should suffice for those more versed in the field and which should also work well for anyone -- new to the area or otherwise-- who is looking to get a basic feel for the theory without getting bogged down in the details. 

As the piece progresses into later sections, the focus of the writing is to convey some of the key ideas and strategies used to prove the results being discussed. For example, the subsections in Section \ref{geometry from combinatorics} are organized by result. Each begins with a summary of the result itself, and then proceeds into a \textit{geometric intuition} subsection whose goal is to give the reader a basic feel for why the result is true. We mostly avoid formal proofs in order to focus on conveying key ideas and intuition.

While we try to give ample references throughout, we stress that this is not the place to find an extensive literature review, and in fact there are many, \textit{many} wonderful papers in this area that have not been mentioned here. Rather, the hope is that the reader will walk away with an increased appreciation for the beautiful interplay between combinatorics and geometry that underlies the study of pseudo-Anosov maps and hyperbolic 3-manifolds.

\subsection{Acknowledgements.} The author thanks Marissa Loving and Brandis Whitfield for carefully reading through a draft of the paper and offering helpful feedback. He also thanks Athanase Papadopoulos for sharing his extensive mathematical-historical knowledge, after which the author was able to write a much more accurate version of footnote 2.

\section{Preliminaries} \label{prelim}

With the exception of the very first subsection below, each subsequent subsection begins with a very brief ``top-level summary'', meant either for experts, or for readers who are not concerned as much with the nitty gritty details. For experts and non-experts alike, the hope is that at least by the standards of a first read, most of the content in Sections \ref{geometry from combinatorics}, \ref{infinite type}, and \ref{Fixed points} can be absorbed only after reading the top-level summaries in this section (and then going back for the more formal details as necessary). 

\subsection{Quasi-comparisons and quasi-isometries}

Given two functions $f, g: A \rightarrow \mathbb{R}$ on a set $A$, we say that $f$ and $g$ are \textit{quasi-comparable} if there exists some $N \in \mathbb{N}$ so that 
\[\forall v \in A,  \frac{1}{N} f(v) - N \leq g(v) \leq N \cdot f(v) + N. \]
Whenever this occurs, we will write $f \sim_{N} g$, where the subscript keeps track of the magnitude of the multiplicative and additive constants in the above inequalities (we will sometimes drop the subscript $N$ for notational convenience). 

Given metric spaces $(X, d_{X}), (Y, d_{Y})$ a map $\phi: X \rightarrow Y$ is a \textit{quasi-isometric embedding}\index{quasi-isometry} if the functions $d_{Y}(\phi( \cdot ), \phi(\cdot)), d_{X}(\cdot, \cdot)$ are quasi-comparable. The map $\phi$ is called a \textit{quasi-isometry} if, in addition, there is some $D \in \mathbb{N}$ so that for any $y \in Y$, there is some $y'$ with $d(y', y) < D$ and so that $y' \in Im(\phi)$. 

If merely $d_{Y}(\phi(x_{1}), \phi(x_{2})) \leq N \cdot d_{X}(x_{1}, x_{2}) + N$ for all $x_{1}, x_{2} \in X$, then $\phi$ is said to be \textit{coarsely Lipschitz}\index{coarse Lipschitz}. 

In section \ref{Fixed points}, we will come across expressions of the form 
\[f \sim_{N} \sum_{Y} [x_{Y}],  \]
where $x_{Y}$ is some constant depending on the summand $Y$, and where $[x_{Y}]$ denotes the quantity 
\[ [x_{Y}] =  \begin{cases} 
      x & x \ge N \\
     0 & x < N.
   \end{cases}
\]

\subsection{Arcs, curves, and surfaces} \label{subsec:acs}

\subsubsection{Top level summary:} \textit{Our surfaces will all be orientable. We are generally interested in homotopy classes of simple closed curves or arcs. In practice, it is common to conflate a curve or an arc with its homotopy class. }

\subsubsection{The details.} Let $S$ be a compact orientable topological $2$-manifold, potentially with boundary. Up to homeomorphism, $S$ is characterized by two invariants: the genus $g$ of $S$, and the number of boundary components $b$. We therefore refer to such manifolds with the notation $S_{g,b}$. If we replace the compactness requirement with the more general assumption of finitely generated fundamental group, such manifolds are fully characterized by adding a third invariant, $p$, the cardinality of a finite set of points removed from some $S_{g,b}$. Hence, by $S_{g,b,p}$ we will mean the unique (up to homeomorphism) orientable surface of genus $g$, $b$ boundary components, and $p$ punctures\index{surface}. 

A \textit{curve}\index{curve} on $S_{g,b,p}$ is a continuous map $\gamma: S^{1} \rightarrow S$; it is \textit{simple}\index{simple!curve} if it is injective, it is \textit{essential} if it is homotopically non-trivial and not homotopic into a tubular neighborhood of a boundary component or puncture, and it is \textit{proper} if its image lies in the interior of the surface. We will focus almost entirely on curves that are simple, proper, and essential. Therefore, unless we explicitly state otherwise, in all that follows, by a \textit{curve} we will mean a curve that satisfies these additional properties. By a \textit{multi-curve}\index{multi-curve}, we will mean a collection of pairwise disjoint and pairwise non-homotopic curves.  

An \textit{arc}\index{arc} is either: 
\begin{enumerate}
    \item a continuous map $\rho: [0,1] \rightarrow S$;
    \item a continuous map $\rho: (0,1) \rightarrow S$ such that $\lim_{t \rightarrow 0^{+}} \rho$ and $\lim_{t \rightarrow 1^{-}} \rho$ are both punctures. 
\end{enumerate}

An arc is said to be \textit{proper} if it is of the second type, or if it is of the first type and $\partial S \cap \rho([0,1]) = \left\{\rho(0), \rho(1) \right\}$. A proper arc is \textit{essential} if it is homotopically non-trivial rel $\partial S$. An arc is \textit{simple} when injective, and \textit{boundary parallel} when it is homotopic into any annular neighborhood of $\partial S$. We will exclusively be concerned with arcs that are proper, essential, and not boundary parallel. Thus, as in the previous paragraph with curves, in all that follows, by an \textit{arc} we will mean an arc that satisfies these additional properties. 

Given $\alpha$ and $\beta$, both of which is either a curve or an arc, the \textit{geometric intersection number}\index{geometric intersection number} $i(\alpha, \beta)$ is defined to be the minimum number of intersections between any curve (or potentially, arc) homotopic to $\alpha$ and any curve (or potentially, arc) homotopic to $\beta$: 

\[ i(\alpha, \beta): = \min \left\{ | \alpha' \cap \beta'|: \alpha' \sim \alpha, \beta' \sim \beta \right\}. \]

We require homotopies of arcs to be homotopies of pairs $([0,1], \partial[0,1]) \hookrightarrow (S, \partial S)$; in particular, homotopies can slide the endpoints of arcs along boundary components.  

The \textit{mapping class group}\index{mapping class group} of $S$, denoted $\mbox{Mod}(S)$, is the group of isotopy classes of orientation preserving homeomorphisms of $S$. When $S$ has boundary components, mapping classes are required to fix boundaries pointwise.  

\subsection{Pseudo-Anosov homeomorphisms} \label{subsec:pA}

\subsubsection{Top-level summary: A pseudo-Anosov map on a surface $S$ with negative Euler characteristic is the type of homeomorphism that most closely resembles the behavior of an Anosov diffeomorphism on a torus: it stretches in one direction and contracts in a transverse direction.}

\subsubsection{The details.} The Nielsen--Thurston classification\index{Nielsen-Thurston} dictates the following trichotomy for homeomorphisms $f: S \rightarrow S$ on a surface $S$ with $\chi(S)<0$ (we refer the reader to \cite{FM} for more details):
\begin{enumerate}
\item $f$ is homotopic to a map of finite order (and in this case, the order is bounded above explicitly in terms of $\chi$);
\item $f$ is not homotopic to a map of finite order, and there exists a multi-curve $\gamma$ on $S$ so that some finite power of $f$ is homotopic to a map that preserves $\gamma$ up to isotopy (and in this case, the absolute value of this power is bounded above explicitly in terms of $\chi$).
\item There exists a singular flat metric $\sigma_{f}$ on $S$ with vertical and horizontal foliations $\mathcal{F}_{v}, \mathcal{F}_{h}$ so that $f$ is homotopic to a map $f'$ which, away from the shared singularities of $\mathcal{F}_{v}, \mathcal{F}_{h}$, stretches along $\mathcal{F}_{v}$ by a factor of some $\lambda >1$ and contracts along $\mathcal{F}_{h}$ by a factor of $1/\lambda$ (equivalently, around a non-singular point, $Df'$ is expressible in local coordinates as a diagonal matrix with diagonal entries $\lambda, 1/\lambda$)\footnote{This characterization of a pseudo-Anosov diffeomorphism is due to Thurston and not to Nielsen, who never reasoned in terms of singular foliations}. 
\end{enumerate}

Those maps in the third category are the so-called \textit{pseudo-Anosov homeomorphisms}. An immediate corollary of this is that a map is pseudo-Anosov precisely when no power of it fixes a homotopy class of simple closed curve. 

Pseudo-Anosovs are by far the most rich and interesting of the three types. For one, they are abundant in the mapping class group: Maher-Tiozzo show that the probability that a random walk on the mapping class group lands on a pseudo-Anosov converges to $1$ exponentially fast in the length of the walk \cite{MaT}, and Choi has shown that pseudo-Anosovs are generic in certain Cayley graphs associated with finite generating sets \cite{Choi}. The constant $\lambda$ is called the \textit{dilatation}\index{dilatation} of the pseudo-Anosov, and it is well-defined in the sense that, at least up to the obvious equivalences, there is only one way to realize a given pseudo-Anosov homotopy class as a map on a singular flat surface as described above.

\subsection{Combinatorial complexes} \label{subsec:CC}

\subsubsection{Top-level summary:} \textit{The curve complex is the simplicial complex\index{simplicial complex} whose simplices correspond to multi-curves on $S$. The simplices of the arc complex correspond to multi-arcs. The vertices of the pants graph correspond to pants decompositions and edges correspond to so-called ``elementary moves.''}

\subsubsection{The details.} The most central character in this story is the so-called \textit{complex of curves}\index{curve complex} of $S_{g,b,p}$, denoted $\mathcal{C}(S)$. It is a simplicial complex whose $0$-simplices correspond to homotopy classes of essential simple closed curves on $S$, and whose $k$-simplices correspond to multi-curves with $k+1$ components. In the language introduced above, a $k$-simplex corresponds to $k+1$ homotopy classes $[\alpha_0],..., [\alpha_k]$ so that 
\[ i(\alpha_n, \alpha_m) = 0, \forall n,m.\]
The curve complex is flag, in the sense that a $k$-simplex appears within it if and only if its $1$-skeleton is already present in the \textit{curve graph} $\mathcal{C}^{1}(S)$, the $1$-skeleton of $\mathcal{C}(S)$.

\begin{figure}
\centering
\vspace{1cm}
\includegraphics[width=\linewidth]{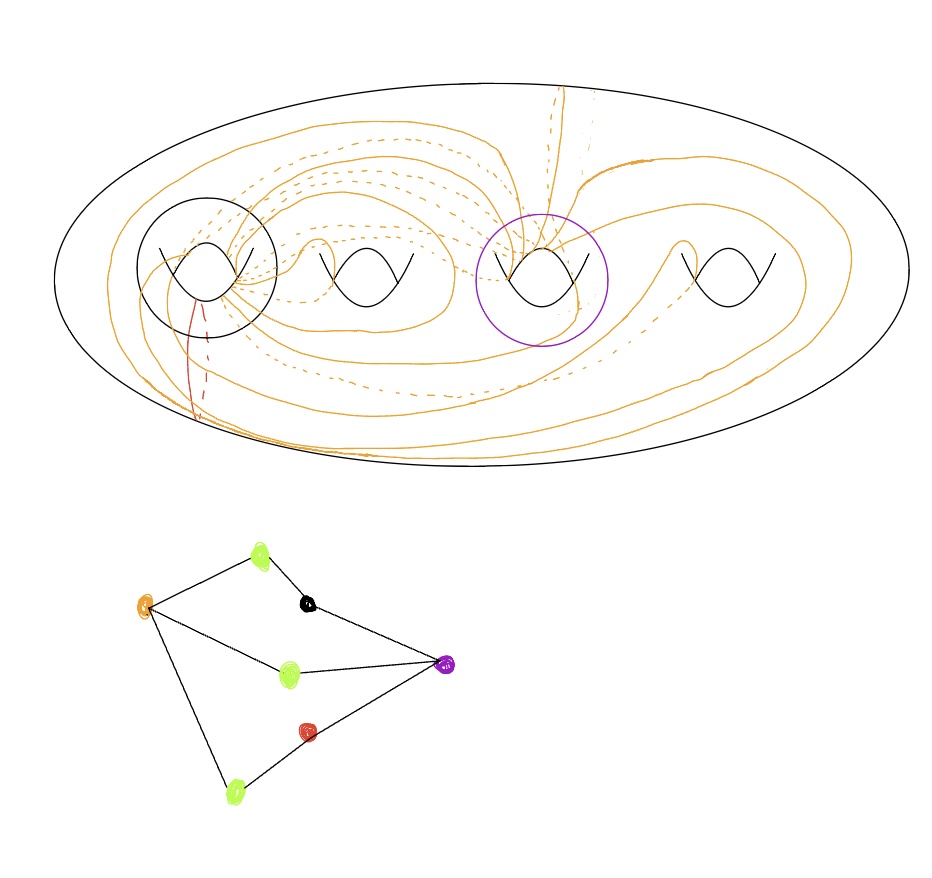}
\caption{Four curves on a genus $4$ surface and the corresponding vertices in the curve complex. We leave it as an exercise for the reader that the orange curve is indeed at a distance of $2$ to each of the black, red, and purple curves (implying the existence of the pictured green vertices).}
\label{fig:curve complex}
\end{figure}

The curve complex was introduced by Harvey \cite{Ha}, and motivated by the various analogies between the mapping class group and $\mathrm{SL}(n, \mathbb{Z})$; through this lens, $\mathcal{C}(S)$ was proposed as a sort of spherical building at infinity for the mapping class group. It is therefore not surprising that its properties encode incredibly deep information about $\mbox{Mod}(S)$. For example, $\mathcal{C}(S)$ has been used to: 
\begin{itemize}
\item solve the conjugacy problem for $\mbox{Mod}(S)$ (\cite{Tao}, \cite{MasurMinskyI});
\item prove a strong Tits alternative for subgroups (\cite{McC}); 
\item bound the asymptotic dimension\index{asymptotic dimension} (\cite{BBF});
\item bound the dimension of a quasi-isometrically embedded copy of $\mathbb{R}^{n}$ and proving that all such \textit{quasi-flats} correspond to free abelian subgroups (\cite{BM}).
\end{itemize}

Note that $\mathcal{C}(S)$ admits an action of $\mbox{Mod}(S)$ by simplicial automorphisms: one simply extends the action on homotopy classes of curves to the higher dimensional skeleta and uses the fact that homeomorphisms preserve disjointness of curves. In addition, the curve complex $\mathcal{C}(S)$ admits a path metric by identifying each simplex with a standard Euclidean simplex with unit length sides. Many of the proofs of the sort of group theoretic results mentioned above boil down to proving this action-- in spite of not being properly discontinuous-- has nice geometric features. This program begins with the foundational result of Masur--Minsky that $\mathcal{C}(S)$ is Gromov hyperbolic\index{Gromov hyperbolic}: for some $\delta>0$, any side of a geodesic triangle is contained in the $\delta$-neighborhood of the union of the other two sides \cite{MasurMinskyI}. 

Before moving on, we will briefly introduce two other related simplicial complexes.

\begin{figure}
\centering
\vspace{1cm}
\includegraphics[scale = 0.2]{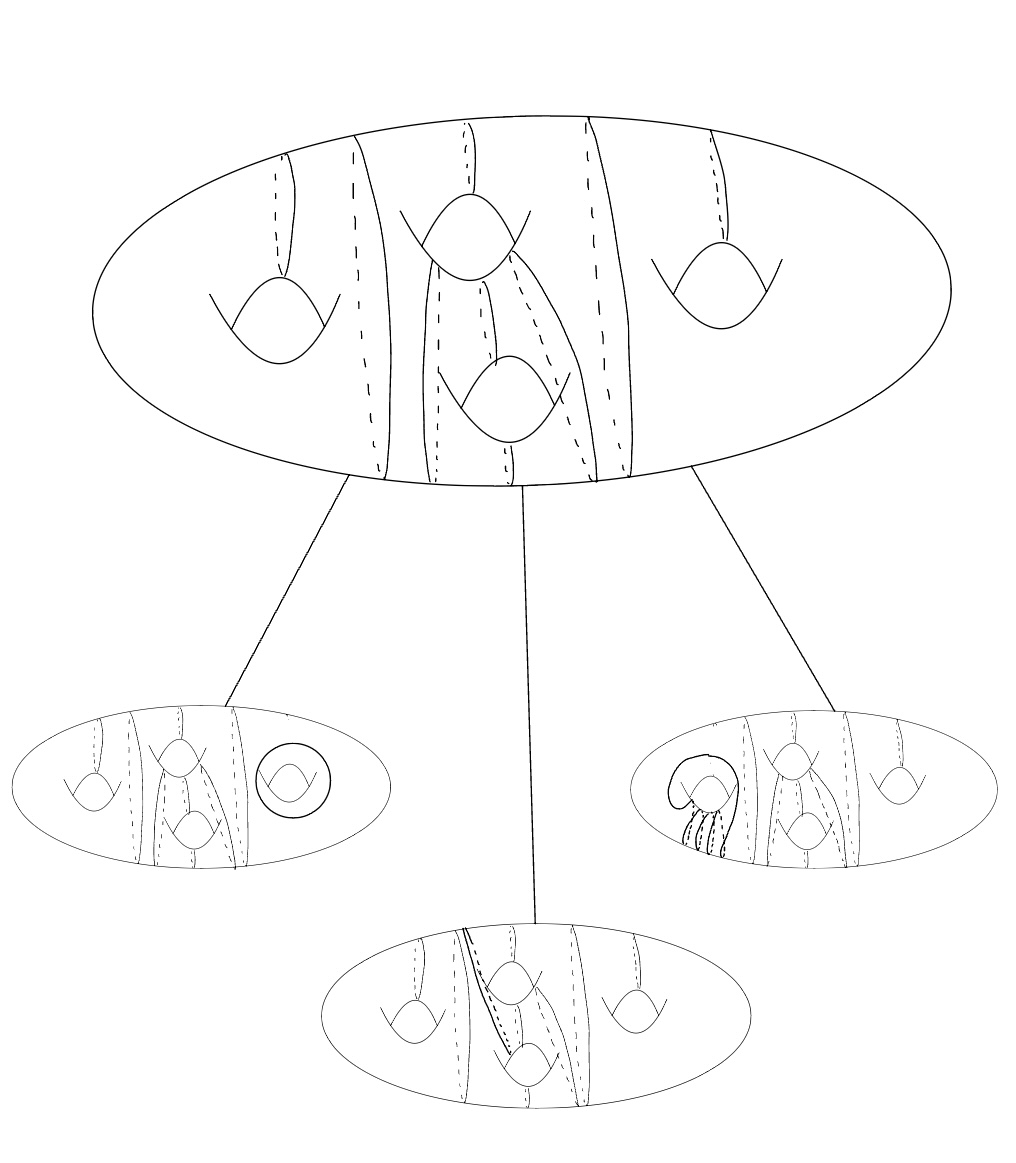}
\caption{Several elementary moves are pictured.}
\label{fig:pants}
\end{figure}

A \textit{pants decomposition}\index{pants decomposition} of $S$, denoted $\mathcal{P}(S)$, is the $0$-skeleton of a maximum--dimensional simplex of $\mathcal{C}(S)$. Equivalently, it is a maximal collection of pairwise disjoint curves on $S$. The \textit{pants graph} is a simplicial graph whose vertices correspond to homotopy classes of pants decompositions, and whose edges correspond to so-called \textit{elementary moves} between decompositions, described as follows: a pants decomposition $P$ is connected to another decomposition $P'$ whenever $P'$ can be obtained from $P$ by deleting a single curve $\gamma$ and replacing it with some curve $\gamma'$ such that $\gamma'$ intersects $\gamma$ minimally amongst all possible replacement curves-- see Figure \ref{fig:pants}.

It is easy to see that in any elementary move as described above $\gamma'$ intersects $\gamma$ either once or twice, depending on whether the connected component of $S \setminus (\mathcal{P} \setminus \gamma)$ which is not of the form $S_{0,p,b}$ where $p+b = 3$, is of the form $S_{1,p,b}$ with $p + b= 1$ or $S_{0,p,b}$ with $p+b = 4$, respectively.

The \textit{arc complex}\index{arc complex} of $S$, denoted $\mathcal{A}(S)$, is defined analogously to $\mathcal{C}(S)$ but using arcs in place of curves for the $0$-simplices: in the event that $S$ has only punctures and no boundary components, higher dimensional simplices correspond to collections of isotopy classes of arcs that can be realized simultaneously disjointly on $S$. We do however need to be careful about conventions when the surface has boundary. Different authors treat this in slightly different ways; in what follows, in accordance with the conventions outlined in Subsection \ref{subsec:acs} we will declare that when $S$ has boundary, isotopies between arcs are allowed to slide endpoints along boundary components.  

One can also consider a complex formed from both arcs and curves, appropriately called the \textit{arc and curve complex} and denoted by $\mathcal{AC}(S)$, where again adjacency corresponds to disjointness on the surface. It is worth pointing out that $\mathcal{C}(S)$ is quasi-isometric to $\mathcal{AC}(S)$ via the map sending each curve to itself. Moreover, there is a coarsely Lipschitz map $\phi: \mathcal{A}(S) \rightarrow \mathcal{AC}(S)$ which comes from sending each arc to itself. It follows from the basic definitions that Gromov hyperbolicity of $\mathcal{A}(S)$ implies the same for $\mathcal{C}(S)$. 

Two decades after Masur--Minsky's original proof of the hyperbolicity of $\mathcal{C}(S)$, Hensel--Przytycki--Webb used this precise line of reasoning to establish a much shorter and far more combinatorial proof of hyperbolicity, by first proving it for $\mathcal{A}(S)$ \cite{HPW}. This new proof established the additional fact that, so long as we focus on the curve \textit{graph} and not the full complex, the $\delta$-constant for hyperbolicity can be taken to be independent of the surface $S$ (this consequence was also deduced through different methods around the same time by the author \cite{Aougab}, Bowditch \cite{B}, and Clay--Rafi--Schleimer \cite{CRS}).

\subsection{Hyperbolic geometry} 

\subsubsection{Top-level summary: Short curves admit long tubular neighborhoods on a hyperbolic surface (the collar lemma). There is a universal constant depending only on the topology of $S$ that controls the minimum length of a pants decomposition on $S$. Geodesic laminations are the objects that can be obtained by taking limits of simple closed geodesics.}

\subsubsection{The details.} We refer the reader to Buser's excellent book for an in--depth discussion of hyperbolic geometry (and basic Teichm{\"u}ller theory) \cite{Bu}. The uniformization theorem\index{uniformization theorem} implies that so long as $g \ge 1$ and $b+p \ge 1$, or $g =0$ and $b+p \ge 3$, $S$ admits a complete Riemannian metric of constant sectional curvature equal to $-1$. When $S$ is of finite type, this metric is of finite area and the Gauss--Bonnet theorem implies that the area depends only on $g,b$, and $p$. Such a metric is said to be \textit{hyperbolic}, owing to the fact that it is modeled locally on the geometry of the hyperbolic plane, $\mathbb{H}^{2}$. Henceforth, when a surface admits a hyperbolic metric we will call it \textit{hyperbolizable}. 

In a hyperbolic metric, every essential homotopy class admits a unique geodesic representative, and geodesic representatives are automatically in minimal position. For these reasons, hyperbolic metrics are very well suited for studying the combinatorics of curves on surfaces. 

A key property of a hyperbolic metric on $S$ is implied by the so-called \textit{collar lemma}\index{collar lemma}: a simple closed geodesic $\alpha$ admits an embedded tubular neighborhood $N(\alpha)$ whose diameter is roughly logarithmic in the reciprocal of the length of $\alpha$. Thus, short curves have very long tubular neighborhoods (this will be made more precise in Subsection \ref{subsec:Margulis} below). 

The \textit{Bers constant}\index{Bers constant} $B= B(S)$ is an upper bound for the length of a shortest pants decomposition for any hyperbolic metric on $S$. Its existence follows from a standard induction argument based on the foundation of the Gauss-Bonnet theorem. One first imagines expanding a round hyperbolic disk on the surface until the moment it achieves a point of self-tangency. Since the total area of $S$ is $2\pi |\chi(S)|$, this occurs when the radius is at most some explicit function of $\chi$ (which grows logarithmically in $|\chi|$ because the area of a hyperbolic disk depends exponentially on its radius). The point of self-tangency signals the existence of an essential simple closed curve, $\alpha_{1}$, whose length is at most the diameter of the disk. Upon finding $\alpha_{1}$, we can cut along it and induct on the topology of the underlying surface to deduce the existence of $B$.

A \textit{geodesic lamination}\index{lamination} on a hyperbolic surface $X$ is a closed subset of $X$ foliated by geodesics. Examples include simple closed geodesics, collections of pairwise disjoint simple closed geodesics, and the Gromov--Hausdorff limit of the orbit of any closed geodesic under successive powers of a pseudo-Anosov diffeomorphism. One can also refer to a subset of a surface (without any metric, hyperbolic or otherwise) as a lamination, exactly when there is a homeomorphism of that surface to a hyperbolic surface which sends the subset to a geodesic lamination.

\subsection{Subsurface projection}

\subsubsection{Top-level summary: The subsurface projection of a curve or an arc is more or less simply the intersection of that curve or arc with the subsurface in question. Care has to be taken to formalize this, especially when it comes to annular subsurfaces since the intersection of a curve with an annulus depends sensitively on where one declares the boundary of that annulus to be. The standard way of dealing with this is to impose a hyperbolic metric, lift to the cover corresponding to the annular subsurface, and then use the natural boundary at infinity of that flaring annular cover. To obtain a unified approach, we will formalize subsurface projection by using this method for all subsurfaces. But the key take-away is that it essentially sends a curve or arc in the full surface to its intersection with a subsurface.}

\subsubsection{The details.} If the curve complex were locally compact and the action of $\mbox{Mod}(S)$ upon it were properly discontinuous, through a standard application of the Schwarz--Milnor lemma one could deduce Gromov hyperbolicity of the mapping class group itself. However, Dehn twists about disjoint curves generate free abelian subgroups and no hyperbolic group contains copies of $\mathbb{Z}^{n}$. Thus, the best one could hope for is a weaker form of hyperbolicity, and indeed, Masur--Minsky use the geometry of $\mathcal{C}(S)$ to initiate a wide-ranging study of the hyperbolic behavior of the large-scale geometry of $\mbox{Mod}(S)$ through so-called \textit{hierarchy paths}\index{hierarchies} \cite{MasurMinskyII}. These paths are built out of geodesics in the curve complexes of various subsurfaces of the full surface $S$-- see subsection \ref{Systole finite type} and Figure \ref{hierarchy} below for more details. The structure enjoyed by the mapping class group as a result of the existence of hierarchy paths has been extensively leveraged for proving all sorts of results (for instance those listed in Subsection \ref{subsec:CC}). And more recently, Behrstock--Hagen--Sisto axiomatized Masur--Minsky's hierarchy machinery and introduced the concept of \textit{hierarchically hyperbolic spaces} (\cite{BSH1}, \cite{BSH2}); this class of course includes mapping class groups, but in addition contains right-angled Artin groups, relatively hyperbolic groups, groups acting on CAT(0) cube complexes, and fundamental groups of many $3$-manifolds.  

While we will not have need to get into the details of the hierarchy machinery (although we will cover \textit{some} of these details in Subsection \ref{Systole finite type}), we will need one of the key tools upon which it is based: so-called \textit{subsurface projection}\index{subsurface projection}. Let $Y \subseteq S$ be an \textit{essential subsurface}, which consists of the following data: 
\begin{itemize}
\item A compact surface $Y$ that is not of the form $S_{b,p}$ with $b+p \leq 1$ or with $b+p =3$. 
\item A continuous map $\iota: Y \rightarrow S$ that is injective and $\pi_{1}$-injective on the interior of $Y$ and which maps each boundary component of $Y$ to an essential curve in $S$. 
\end{itemize}
To each such subsurface $Y$, there corresponds a map 
\[ \pi_{Y}: \mathcal{C}^{0}(S) \rightarrow 2^{\mathcal{AC}^{0}(Y)}\]
from the set of essential curves in $S$ to the power set of the $0$-skeleton of the arc and curve graph of $Y$, defined as follows. 

\begin{figure}
\centering
\vspace{1cm}
\includegraphics[scale=0.2]{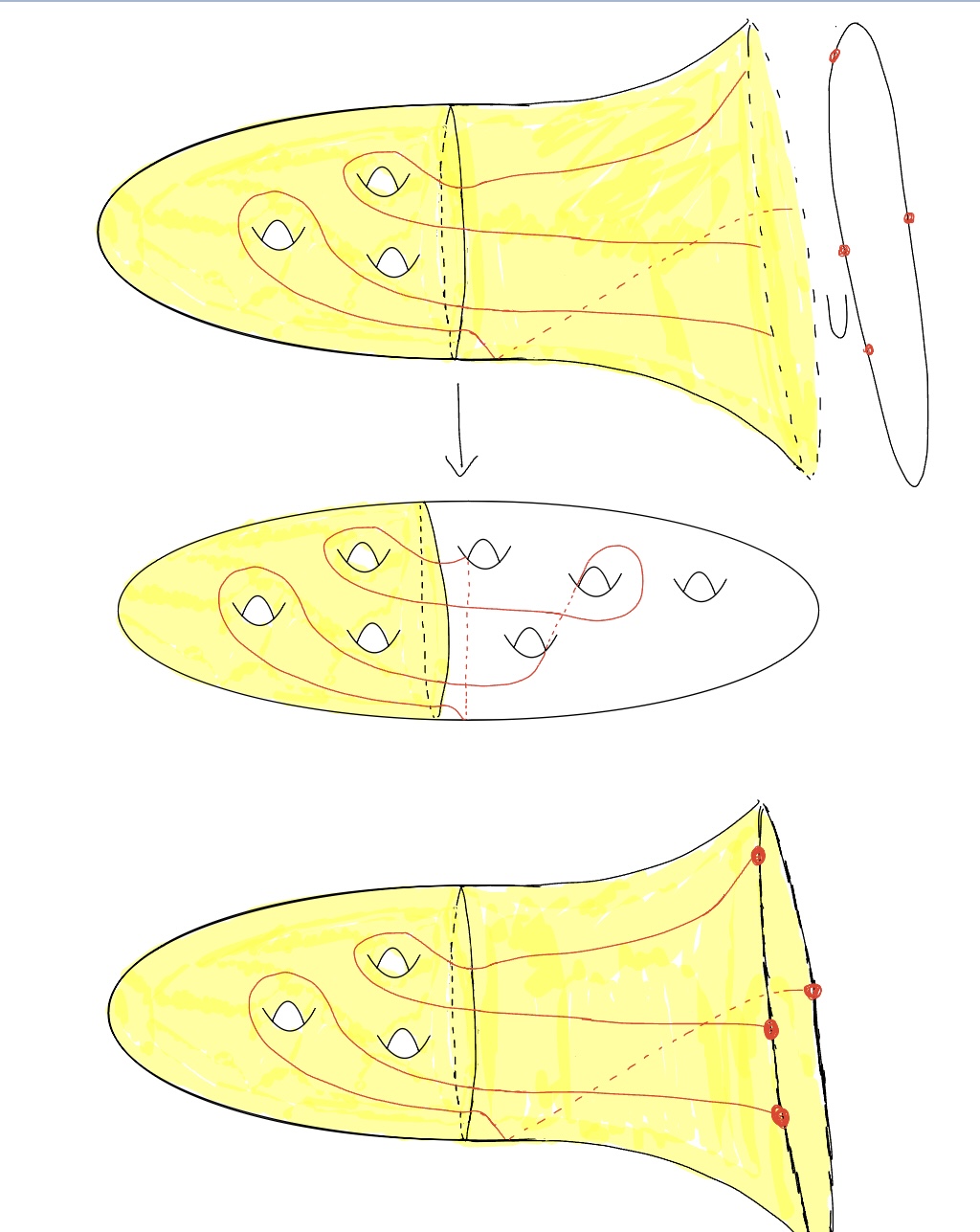}
\caption{A curve together with its subsurface projection to a yellow subsurface. The top figure represents the cover associated with this subsurface, which can be compactified using the boundary at infinity in the universal cover, yielding the bottom picture.}
\label{fig:subsurface}
\end{figure}

Fix once and for all a complete and finite area hyperbolic metric on $S$. Then associated with $Y$ is a covering space $p_{Y}: \widetilde{Y} \rightarrow S$. The universal covering map from $\mathbb{H}^{2}$ to $\widetilde{Y}$ extends continuously to the boundary at infinity, yielding a compactification of $\widetilde{Y}$ to a surface $\overline{Y}$ with boundary. We define $\pi_{Y}(\alpha)$ to be the closure of $p_{Y}^{-1}(\alpha) \subset \overline{Y}$- see Figure \ref{fig:subsurface}.

One can often interpret $\pi_{Y}$ in a coarser and more combinatorial sense as simply intersecting a representative of $\alpha$ (that is in minimal position with $\partial Y$) with $Y$. This approach fails however when $Y$ is an annulus since in that case, the two boundary components of $Y$ are isotopic to each other and so at least one of them is necessarily not geodesic, so there is no canonical representative for $Y$.

The key observation that, in practice, allows one to reinterpret $\pi_{Y}$ as a map between full complexes and not merely a map from vertices to a power set, is that the diameter of $\pi_{Y}(\alpha)$ in $\mathcal{AC}(Y)$ is uniformly bounded. Indeed, the components of $\iota^{-1}(a)$ are all disjoint from one another by virtue of $\alpha$ being simple, and the various representatives in $A_{\alpha}$ can be assumed to be disjoint from one another within the interior of $\iota(Y)$. For the same reason, if one starts with a \textit{pair} of disjoint curves $\alpha_{1}, \alpha_{2}$, the diameter of the union of their images under $\pi_{Y}$ will also be uniformly bounded. We summarize this in the following lemma, due originally to Masur-Minsky \cite{MasurMinskyII}: 

\begin{lemma} \label{lem:lipchitz projections} The map $\pi_{Y}$ induces a coarsely Lipshitz map (which by a slight abuse of notation we refer to as $\pi_{Y}$ as well) 
\[ \pi_{Y}: \mathcal{C}(S) \rightarrow \mathcal{AC}(Y). \]
\end{lemma}

Given curves $\alpha_{1}, \alpha_{2}$ and an essential subsurface $Y$, by $d_{Y}(\alpha_{1}, \alpha_{2})$ we will mean the diameter in $\mathcal{AC}(Y)$ of $\pi_{Y}(\alpha_1) \cup \pi_{Y}(\alpha_{2})$.

While we omit the formal details, if $\lambda$ is a lamination on $S$ we can also define $\pi_{Y}(\lambda)$
in a similar fashion: intersect $\lambda$ with $Y$ and keep track only of the resulting pairwise disjoint homotopy classes of arcs and curves that arise.

\subsection{Teichm{\"u}ller space}

\subsubsection{Top-level summary: The Teichm{\"u}ller space parameterizes marked complete hyperbolic metrics on a fixed surface $S$. It can be topologized using either the Teichm{\"u}ller metric (in which the distance between surfaces is proportional to the logarithm of the optimal dilatation for any (marked) quasi-conformal homeomorphism between them) or the Weil-Petersson metric (in which distance is given by means of a Riemannian inner product coming from integrating the product of a pair of quadratic differentials against the hyperbolic area form.)}

\subsubsection{The details.} The \textit{Teichm{\"u}ller space}\index{Teichm{\"u}ller space} $\mathcal{T}(S)$ of a finite type surface $S$, is, as a set, a quotient $\left\{ ([\phi], \sigma)\right\}/\sim$ of the set of pairs $([\phi], \sigma)$ where: 
\begin{itemize}
\item $\sigma$ is a $2$-manifold homeomorphic to $S$ equipped with a hyperbolic metric, and  
\item $\phi: S \rightarrow \sigma$ is a homotopy class of homeomorphisms.
\end{itemize}
The equivalence relation $\sim$ is given by declaring that pairs $([\phi_1], \sigma_1), ([\phi_2], \sigma_2)$ are equivalent precisely when there exists an isometry $j: \sigma_1 \rightarrow \sigma_2$ and representatives $\phi_1 \in [\phi_1], \phi_2 \in [\phi_2]$ so that $j$ is homotopic to $\phi_1 \circ \phi_2^{-1}$.

A consequence of the uniformization theorem is that $\mathcal{T}(S)$ can also be described as a set of pairs where the second factor is a Riemann surface homeomorphic to $S$, and where the equivalence is defined by the existence of a conformal automorphism in the appropriate homotopy class.  

The Teichm{\"u}ller space can be topologized in a number of different but equivalent ways, and when one equips $\mathcal{T}(S)$ with this topology, it becomes homeomorphic to a ball of dimension $6g-6 + 2(b+p)$. We briefly describe two metrics on $\mathcal{T}(S)$ that give rise to this topology: 

\begin{enumerate}
\item \textbf{The Teichm{\"u}ller metric:} Given two Riemann surfaces $R_{1}, R_{2}$ homeomorphic to $S$ and a diffeomorphism $h: R_{1} \rightarrow R_{2}$, one can consider the \textit{quasi-conformal dilatation}, $\kappa_{h}$ of $h$, which, loosely speaking, considers the maximum amount by which an infinitesimal circle in $TR_{1}$ is warped (this can be quantified by measuring the eccentricity of the ellipse in the image of such a circle -- we refer the reader to Ahlfors for a formal discussion \cite{Ahl}). When this quantity is finite, the map is said to be \textit{quasi-conformal}. Since conformal maps send circles to circles, the deviation of the dilatation from $1$ measures the extent to which a quasi-conformal map is not conformal. 

Teichm{\"u}ller's\footnote{We make the brief historical note that Teichm{\"u}ller was a virulent Nazi. For a time, there was a sort of debate in the community about whether to rename the various theorems and spaces to which we refer as ``Teichm{\"u}ller'' in order not to bestow any honor on the man, whose memory clearly deserves none. Some have argued that engaging in such a renaming would also more accurately reflect the mathematical historical reality by honoring more significant contributions due to, for instance, Fricke and Klein. While the author is partial to some of these arguments, it was brought to his attention by Athanase Papadopoulos that in good faith, one really can not deny the extensive and revolutionary mathematical achievements of Teichm{\"u}ller in this area, for instance: 

He introduced the notion of a marked surface and of the Teichm{\"u}ller space itself, including its topology and complex structure. Teichm{\"u}ller's theorem establishing an extremal quasiconformal homeomorphism between any two marked Riemann surfaces is of course the basis for the Teichm{\"u}ller metric. He was also the first to describe the tangent space via quadratic differentials, and the first to investigate the $SL(2,\mathbb{R})$ orbit of a corresponding tangent vector, an object now known as a Teichm{\"u}ller disk. He reinterpreted the Nielsen realization problem of finding a hyperbolic surface with a given finite isometry group in terms of finding a fixed point for the action of the corresponding finite subgroup of the mapping class group on the Teichm{\"u}ller space, which is of course the interpretation eventually used by Kerckhoff to resolve the problem \cite{Ker}. 

Suffice to say: the argument that Teichm{\"u}ller's name should be erased from mathematical history \textit{on the basis} that his contributions to the theory were outmatched by others, does not hold water. Besides this, the author feels that this debate has happened after the proverbial train has more or less left the station. At this point it would probably cause far more confusion than necessary to actually successfully execute a renaming process, especially to younger mathematicians struggling to learn the theory under a new set of names while everyone their senior still uses the old.

Perhaps the answer is not to erase, but to remember \textit{fully}: Teichm{\"u}ller was both a tremendously influential mathematician and an unforgiveable Nazi, and as such, he is but one of many examples of the reality that mathematical ability does not necessarily correlate with moral or spiritual fortitude.} key theorem is that for any two homeomorphic Riemann surfaces and any specified homotopy class of homeomorphisms $[\phi]$ between them, there exists a unique representative $f_{[\phi]}$ that minimizes the quasi-conformal dilatation. Then one can define a map 
\[ d_{\mathcal{T}} : \mathcal{T}(S) \times \mathcal{T}(S) \rightarrow \mathbb{R}\]
by 
\[ d_{\mathcal{T}}(([\phi_1], \sigma_1), ([\phi_2], \sigma_2)) = \log(\kappa_{f_{\phi_1 \circ \phi_2 ^{-1}}}). \]
This is a metric on $\mathcal{T}(S)$ called the \textit{Teichm{\"u}ller metric}; it is Finsler\index{Finsler} in the sense that it is induced by a smoothly varying norm on tangent spaces. 

\item \textbf{The Weil-Petersson metric:} One can identify the cotangent space to a point $X \in \mathcal{T}(S)$ with the space of \textit{quadratic differentials}\index{quadratic differential} on $X$, objects that are expressible in the form $\phi(z)dz^{2}$ in terms of a local coordinate $z$ and a holomorphic function $\phi$. This identification gives rise to a Hermitian form on $T_{X}\mathcal{T}(S)$: define $\langle v, w \rangle$ to be 
\[ \int_{X} \phi_{v} \overline{\psi_{w}},\]
where $\phi_v(z) dz^{2}$ (resp. $\psi_w (z) dz^{2}$) is the quadratic differential dual to $v$ (resp. $w$) and the integral is taken with respect to the area form coming from the hyperbolic metric on $X$. The associated Riemannian metric on $\mathcal{T}(S)$ is known as the \textit{Weil-Petersson metric}\index{Weil-Petersson}. It is CAT(0) and metrically incomplete; the completion points correspond to \textit{noded hyperbolic surfaces}, which can be obtained as a limit of hyperbolic metrics along which the length of a given multi-curve goes to $0$. A formal discussion of these ideas can be found in the extensive work of Wolpert on this topic (\cite{Wolpert1}, \cite{Wolpert2}). 

\subsection{Dynamical interpretations for the WP metric.} Even though what follows in this subsection has very little to do with the rest of the piece and will not come up again later, the author simply can not resist mentioning an alternative way to understand the Weil--Petersson metric that foregoes any complex analysis. 

Start with some $X \in \mathcal{T}(S)$ and fix some $L>0$. Consider the set $N_{X}(L)$ of all closed geodesics on $X$ with length at most $L$. One can then form the function 
\[ \mathcal{F}_{L}: \mathcal{T}(S) \rightarrow \mathbb{R}, \]
 which computes the \textit{average} length of all geodesics in $N_{X}(L)$ as a function of $X \in \mathcal{T}(S)$. Thurston proved the remarkable theorem that, as $L \rightarrow \infty$, the Hessians of $\mathcal{F}_{L}$ converge (in an appropriate sense) to a positive definite form on $T_{X}\mathcal{T}(S)$, which varies smoothly in $X$. Amazingly, the corresponding Riemannian metric is nothing but a multiple of the Weil--Petersson metric. A proof of this can be found in a paper of Wolpert \cite{Wolpert3}.  

 Therefore, one can assign significant geometric meaning to the Weil--Petersson form: given an infinitesimal direction, it-- in a sense-- measures the second directional derivatives of an \textit{average} length function on $X$. This interpretation has allowed the introduction of Weil--Petersson like metrics on other moduli spaces; the theory of these generalized Weil--Petersson forms is captured by the so-called \textit{thermodynamic formalism}-- we refer the reader to McMullen for more details \cite{McM}.  For example, Policott--Sharp used this approach to introduce a Weil--Petersson like metric on the Culler--Vogtmann outer space \cite{PS}, and the author (jointly with Clay and Rieck) \cite{ACR} showed that it indeed shares many of the global properties of the classical Weil--Petersson metric on $\mathcal{T}(S)$. 

\end{enumerate}

\subsection{Hyperbolic $3$-manifolds} \label{subsec:Margulis}

\subsubsection{Top-level summary: In a hyperbolic $3$-manifold, short closed geodesics admit thick tubular neighborhoods. Many $3$-manifolds admit hyperbolic metrics, and any closed $3$-manifold virtually fibers over the circle. While a hyperbolic metric on a closed $3$-manifold is unique up to isometry (in stark contrast to the $2$-dimensional case), there are uncountably many complete hyperbolic metrics on a manifold homeomorphic to $S \times \mathbb{R}$ where $S$ is some surface itself admitting a hyperbolic metric.}

\subsubsection{The details.}The celebrated Margulis Lemma states that the set of elements in a semisimple Lie group translating a given point in the corresponding homogeneous space a sufficiently small amount, generates a virtually nilpotent subgroup. A consequence of this -- plus some standard structural results of the isometry groups of hyperbolic $n$-space-- imply that there is some constant $\epsilon_{n}$ depending only on the dimension, $n$, such that set of all points in a closed hyperbolic $n$-manifold $M^{n}$ with injectivity radius at most $\epsilon_{n}$ decomposes into a disjoint union of tubular neighborhoods of geodesics. We refer to these as $\textit{Margulis tubes}$, and the constant $\epsilon_{n}$ as the $n$-dimensional \textit{Margulis constant}; in dimension $2$, a Margulis tube has the geometry and topology of a flaring annulus with the geodesic at the center. In dimension $3$, it takes the form of a solid torus with the corresponding geodesic as the longitude. We refer the reader to \cite{Kap} for a more in--depth discussion. 

By $M_{\le \epsilon_{n}}$ we will mean the subset of $M$ consisting of points with injectivity radius at most $\epsilon_{n}$. When $n=2,3$ and $M$ has finite volume, $M_{\le \epsilon_{n}}$ consists of a disjoint union of Margulis tubes and \textit{cusp neighborhoods}, a quotient of a horoball by the action of either a rank $1$ or $2$ abelian subgroup of $\mbox{Isom}^{+}(\mathbb{H}^{n})$. 

Thurston's celebrated hyperbolization theorem implies that any closed Haken\index{Haken} $3$-manifold (one which admits a properly embedded $2$-sided and $\pi_{1}$-injective closed surface) which is atoroidal\index{atoroidal} (meaning that it does not admit a $\pi_{1}$-injective map of a torus) admits a hyperbolic metric. A special case of this corresponds to the mapping torus of a pseudo-Anosov diffeomorphism $\phi$ on a finite type hyperbolizable surface $S$ \cite{Thurston3}. On the other hand, it is a consequence of the virtual fibering theorem that this case is in fact not so ``special'': every closed hyperbolic $3$-manifold is finitely covered by such a fibered manifold \cite{Agol3}.  

Mostow--Prasad Rigidity\index{Mostow-Prasad rigidity} implies that when a $3$-manifold admits a complete finite-volume hyperbolic metric, it is unique up to isometry. Therefore, the hyperbolic metrics described in the previous subsection on the manifolds fibering over the circle are unique, and so geometric properties of these metrics can be associated as invariants of the corresponding monodromy mapping classes. 

On the other hand, if $M$ is for instance homeomorphic to $S \times \mathbb{R}$, there are many pairwise non-isometric hyperbolic metrics with which one can equip $S$. The so-called \textit{convex cocompact}\index{convex cocompact} metrics-- ones which come from the quotient of $\mathbb{H}^{3}$ by a discrete group of isometries with the property that the action restricted to the convex hull of the limit set is cocompact (equivalently, those manifolds admitting a convex sub-manifold whose inclusion induces an isomorphism on the level of fundamental groups) -- are parameterized using Bers' simultaneous uniformization theorem \cite{Bers}. They correspond one-to-one with points in $\mathcal{T}(S) \times \mathcal{T}(\overline{S})$, where $\overline{S}$ denotes $S$ but equipped with the opposite orientation. Loosely speaking, one realizes this correspondence as follows. 

The limit set of the aforementioned discrete group of isometries will be a Jordan curve on $\partial_{\infty} \mathbb{H}^{3} = S^{2}$. It separates the $2$-sphere into two simply connected domains of discontinuity, $\Omega_{+}$ and $\Omega_{-}$. The quotient of the action of the group on either such domain produces a marked Riemann surface, and the content of Bers' theorem is that the pair of such surfaces determines the hyperbolic metric of the entire manifold. 

For an example of a metric on $S \times \mathbb{R}$ which is not convex cocompact, one can consider the infinite cyclic cover of the mapping torus $M_{\phi}$ of a pseudo-Anosov diffeomorphism, $\phi$. Indeed, starting with a closed geodesic $\alpha$, the images under higher and higher powers of $\phi$ produce a sequence of geodesics which exit one of the two ends of $M$. 

Work of Bonahon \cite{Bon} and Canary \cite{Can2} imply that these are essentially the only two possibilities for an end of a hyperbolic $3$-manifold $M$ with finitely generated fundamental group; it is either convex cocompact and is thus foliated by surfaces which, after appropriately rescaling, converge to a hyperbolic metric ``at infinity'', or it is \textit{degenerate} in the sense that there is a sequence of closed geodesics of $M$ that exit the end. In the latter case, these geodesics converge to a lamination on $S$, called the \textit{ending lamination}\index{ending lamination} of that particular end.

\subsection{Pleated and simplicial surfaces} \label{subsec:pleated}

\subsubsection{Top-level summary: A pleated surface is essentially a $\pi_{1}$-injective surface in a hyperbolic $3$-manifold $M$ so that the pull-back of the ambient metric on $M$ induces a complete hyperbolic metric on the surface. The map of the surface into $M$ is ``bent'' along some geodesic lamination. By specifying a pair of geodesics in $M$, one can often ``sweep'' through $M$ by a family of pleated surfaces which start at one geodesic and end at the other.}

\subsubsection{The details.} Given $M$ a hyperbolic manifold homeomorphic to $S \times \mathbb{R}$ and $\lambda \subset S$ a lamination, a \textit{pleated surface}\index{pleated surface} with \textit{pleating lamination} $\lambda$ is a pair $(f, \sigma)$ where: 
\begin{enumerate}
\item $\sigma$ is a hyperbolic surface homeomorphic to $S$;
\item $f: \sigma \rightarrow M$ is a $1$-Lipschitz $\pi_{1}$-injective map sending $\sigma$ to its geodesic representative in $M$. 
\end{enumerate}
Thurston first observed that given any $\lambda$ and any such $M$ admitting a geodesic representative of $\lambda$, there always exists a pleated surface with pleating locus $\lambda$ \cite{Thurston1}. 

Pleated surfaces arise very naturally in the study of hyperbolic $3$-manifolds. For example, if $M$ is quasi-Fuchsian then the two components of the boundary of its convex core are naturally pleated surfaces. The content of Thurston's so-called bending conjecture, proved by Dular-Schlenker \cite{DS}, is that the data of these two pleating loci determine the geometry of the entire manifold $M$. 

One could also begin with a triangulation $T$ of $S$, and seek a hyperbolic metric and a $1$-Lipschitz $\pi_{1}$-injection that sends the edges of $T$ to geodesic segments in $M$. Such an object is called a \textit{simplicial hyperbolic surface}-- see \cite{Canary} for more details. 

If $M \approx S \times \mathbb{R}$ is a hyberbolic $3$-manifold and if $\alpha, \beta$ is a pair of simple closed geodesics, a \textit{sweep-out} in $M$ from $\alpha$ to $\beta$ is a map $H: [0,1] \times S \rightarrow M$ so that
\begin{itemize}
\item For each $t$, $H|_{t}$ is a pleated (or perhaps simplicially hyperbolic) surface; 
\item For $t = 0$ (resp. $t= 1$), $\alpha$ (resp. $\beta$) is in the pleating locus of $H|_{t}$ (in the event that we use the simplicially hyperbolic version of a sweep-out, we require instead that the triangulation at time $0$ intersects $\alpha$ minimally amongst all triangulations with a fixed vertex set, and similarly for $t=1$ and $\beta$). 
\end{itemize}

\section{Geometry from combinatorics} \label{geometry from combinatorics}

In this section, we describe some of the results from the literature that allow one to read off geometric information from the combinatorics of one of the simplicial complexes described in \ref{subsec:CC}.

Let $\phi$ be a pseudo-Anosov diffeomorphism on a finite type orientable surface $S$ without boundary. One can then form the mapping torus\index{mapping torus} $\M_{\phi}$, which topologically is the $3$-manifold 
\[ M_{\phi} = (S \times [0,1])/(x,0) \sim (\phi(x), 1), \]
which, by Thurston's hyperbolization theorem, carries a complete hyperbolic metric. Punctures on $S$ will give rise to rank-$2$ cusps in this $3$-manifold, a neighborhood of which deformation retracts to a torus which is the quotient of a horosphere in $\mathbb{H}^{3}$ by a copy of $\mathbb{Z}^{2}$. 

In each subsection below, a result relating geometry and combinatorics  is first described, and there is then a subsubsection that summarizes the geometric intuition behind a proof of that result. The reader is encouraged to note the unifying themes that bind these various arguments together: pleated and simplicially hyperbolic surfaces, sweep-outs, and the machinery developed by Masur-Minsky for studying the geometry of the curve complex and the mapping class group.

\subsection{Electric circumference} \label{EC} Letting $\epsilon= \epsilon_{3}$ denote the $3$-dimensional Margulis constant and $\gamma \in \pi_{1}(M_{\phi})$, for some fixed $\delta < \epsilon$ let $\mathbb{T}_{\delta}(\gamma)$ denote the Margulis tube about $\gamma$ which is the connected component of the set of points in $M_{\phi}$ which contains $\gamma$ and with injectivity radius less than $\delta$. There are three mutually disjoint possibilities regarding the topology of $\mathbb{T}_{\delta}(\gamma)$: 

\begin{enumerate}
\item It is empty, in the event that every representative of $\gamma$ is sufficiently long (for instance, at least $2 \cdot \delta$).
\item It is homeomorphic to a solid torus whose core curve is the geodesic representative for $\gamma$, in the event that the conjugacy class of $\gamma$ acts loxodromically on $\mathbb{H}^{3}$. 
\item It is homeomorphic to a solid torus with its core drilled out (and therefore it deformation retracts to a torus) in the event that $\gamma$ acts parabolically. 
\end{enumerate}

Let $\rho: [0,1] \rightarrow M_{\phi}$ be a $d$-rectifiable path. There is then a partition $0= t_{0} < t_{1} < t_{2} < ... < t_{n}= 1$ so that either $\rho|_{(t_{i}, t_{i+1})}$ lies outside of all Margulis tubes, or, $\rho|_{(t_{i}, t_{i+1})}$ lies in some $\mathbb{T}_{\delta}(\gamma)$. By choice of $\delta$, distinct $\delta$-Margulis tubes have disjoint boundaries, and therefore the sub-intervals just described must alternate between thick and thin: if $\rho$ lies in a $\delta$-tube between $t_{i}$ and $t_{i+1}$, then it is in the thick part between times $t_{i-1}$ and $t_{i}$, and between $t_{i+1}$ and $t_{i+2}$. 

Letting $\ell$ denote the hyperbolic arc length function on $M_{\phi}$, we will define a new arc length $\ell_{\delta}$-- called the $\delta$-\textit{electrification} of $\ell$-- as follows: 

\[ \ell_{\delta} (\rho) = \sum_{i=0}^{n-1} \ell_{\delta}(\rho|_{(t_{i}, t_{i+1})}), \]
where $\ell_{\delta}(\rho|_{(t_{i}, t_{i+1})}) = \ell(\rho|_{(t_{i}, t_{i+1})})$ when $\rho$ lies outside of all Margulis tubes in the interval $(t_{i}, t_{i+1})$, and where $\ell_{\delta}(\rho|_{(t_{i}, t_{i+1})}) = 1$ otherwise. 

The idea here is that we are effectively coning off every $\delta$-tube in $M_{\phi}$. Letting $d$ denote the hyperbolic metric on $M_{\phi}$, we can then define the $\delta$-\textit{electrification} of $d$, denoted $d_{\delta}$, to be the path metric associated with this arc length function. 

Finally, the $\delta$-\textit{electric circumference} of $M_{\phi}$, denoted $\mbox{circ}_{\delta}(M_{\phi})$, is defined to be 
\[ \inf(\ell_{\delta}(\rho): \rho : S^{1} \rightarrow M \hspace{1 mm} \mbox{and} \hspace{1 mm} [\rho] \notin \ker(\kappa: \pi_{1}(M)\rightarrow \mathbb{Z}) ),  \]
 where $\kappa$ is the homomorphism to $\mathbb{Z}$ induced by the fibration\index{fibration} of $M$ over $S^{1}$. Informally, the electric circumference measures the hyperbolic arc length of a loop traversing the monodromy of $M_{\phi}$, after ignoring the length that takes place in Margulis tubes and keeping track only of the length occurring in the thick part and the number of tubes through which the loop travels. 

The following theorem, which can be attributed to Brock--Bromberg \cite{BB} (or to Brock--Canary--Minsky \cite{BCM} since it follows also from the general model manifold machinery), relates the electric circumference to distance in the curve complex: 

\begin{theorem} \label{thm:electric} There is a constant $K$ depending only on $\chi(S)$ and $\delta$, so that 
\[ \tau_{\mathcal{C}}(\phi) \sim_{K} \textnormal{circ}_{\epsilon}(M_{\phi}). \]
where $\tau_{\mathcal{C}}(\phi)$ denotes the \textnormal{asymptotic translation length} of $\phi$ in the curve complex:
\[\tau_{\mathcal{C}}(\phi) = \lim_{n \rightarrow \infty}\frac{d_{\mathcal{C}}(v, \phi^{n}(v))}{n}, \]
where $v \in \mathcal{C}^{0}(S)$ is chosen arbitrarily. 
\end{theorem}

\subsubsection{Geometric intuition:} Why should one expect the curve complex to capture anything about the electric circumference? To understand the connection, imagine a path through $\mathcal{C}^{1}(S)$ of the form $\left\{v_{0}, v_{1},....,v_{n} = \phi(v_{0}) \right\}$. Since $v_{i}$ is disjoint from $v_{i+1}$, the union $v_{i} \cup v_{i+1}$ is a multi-curve on $S$, and hence a lamination. Thus there is a pleated surface $(f_{i}, \sigma_{i})$ with pleating locus $v_{i} \cup v_{i+1}$. 

\begin{figure}
\centering
\vspace{1cm}
\includegraphics[scale= 0.2]{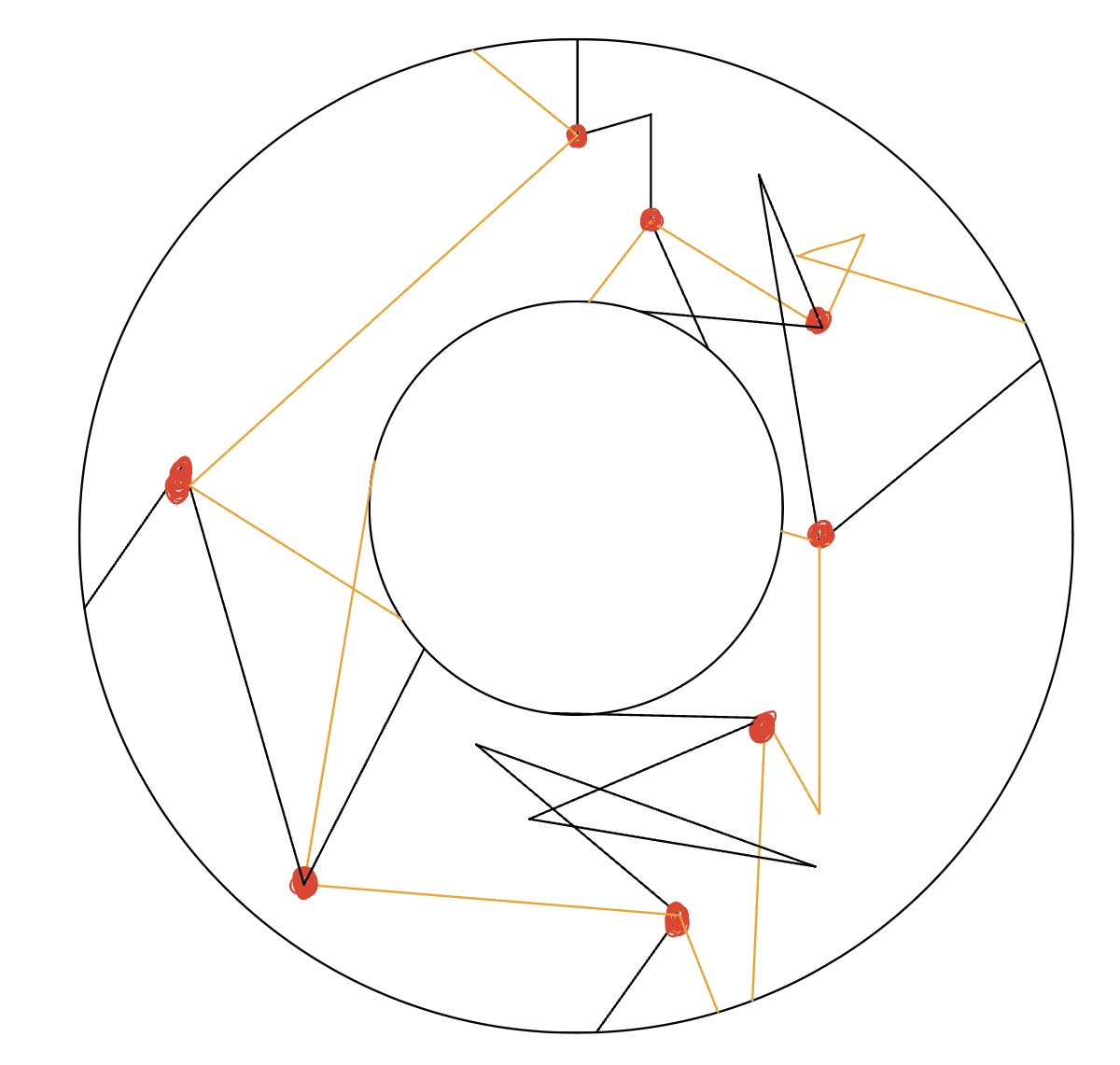}
\caption{Each black or orange piecewise linear curve represents a pleated surface in the mapping torus. The red dots represent simple closed geodesics that lie within the pleating locus of a given pleated surface. Starting with the red geodesic $\alpha$ located highest on the page, there is a black pleated surface containing both it, and a curve $\alpha_{1}$ disjoint from $\alpha$. Then, one discovers a pleated surface (pictured in orange) containing both $\alpha_{1}$ and a curve $\alpha_{2}$ disjoint from $\alpha_{1}$. The colors alternate between black and orange in order to make it easier to see the ``chain''. We also draw some of the pleated surfaces in a way that emphasizes to the reader that the maps need not be embeddings.}
\label{fig:chain}
\end{figure}

One thus obtains a ``chain'' of pleated surfaces\index{pleated surface} in $M$ such that successive surfaces in the chain share a geodesic in common-- see Figure \ref{fig:chain}. The final geometric ingredient to keep in mind is a simple application of the Gauss-Bonnet theorem: the diameter of a hyperbolic surface is bounded solely in terms of its Euler characteristic and its injectivity radius\footnote{To understand this, imagine that the diameter of a hyperbolic surface with large injectivity radius, is also large. There is then a very long geodesic segment in the surface achieving the minimum length between its endpoints. For each point along this path, imagine the shortest closed loop intersecting that point. These loops sweep out a cylinder of large area, equal to roughly the length of the path times the injectivity radius}. Since the map from each pleated surface into $M$ is $1$-Lipschitz, points of small injectivity are mapped to points of small injectivity in $M$; therefore, the $\delta$-electric distance between the geodesic representatives $v_{i}$ and $v_{i+1}$ is uniformly bounded from above in terms only of $\chi(S)$. The upshot is that, because $\phi(v_{0}) = v_{n}$, there is a loop in $M$ not in the kernel of the map from $\pi_{1}(M)$ to $\mathbb{Z}$, whose $\delta$-electric arclength is proportional to $n$.

The converse direction-- that is, bounding the electric circumference from below in terms of $\tau_{\mathcal{C}}(\phi)$, is more complicated and originally relied on compactness and limiting arguments which are considerably harder to make concrete, but we can still give a vague sketch of the idea. The key tool is the notion of a sweep-out (see subsection \ref{subsec:pleated}). Working in the infinite cyclic cover $\widetilde{M_{\phi}}$ and letting $\alpha, \beta$ be an arbitrary pair of simple closed geodesics, one imagines a sweep-out $H: S \times [0,1] \rightarrow \widetilde{M_{\phi}}$ from $\alpha$ to $\beta$. 
 In the context of electric circumference, we first consider an electric geodesic: an arc $\omega: [0,1] \rightarrow \widetilde{M_{\phi}}$ starting on a lift of the fiber $S$ and ending on $\phi(S)$ whose electric arclength is minimal. Consider then a sweep-out where $\beta = \phi(\alpha)$. 

 \vspace{2 mm}

\begin{figure}
\centering
\vspace{1cm}
\includegraphics[scale = 0.2]{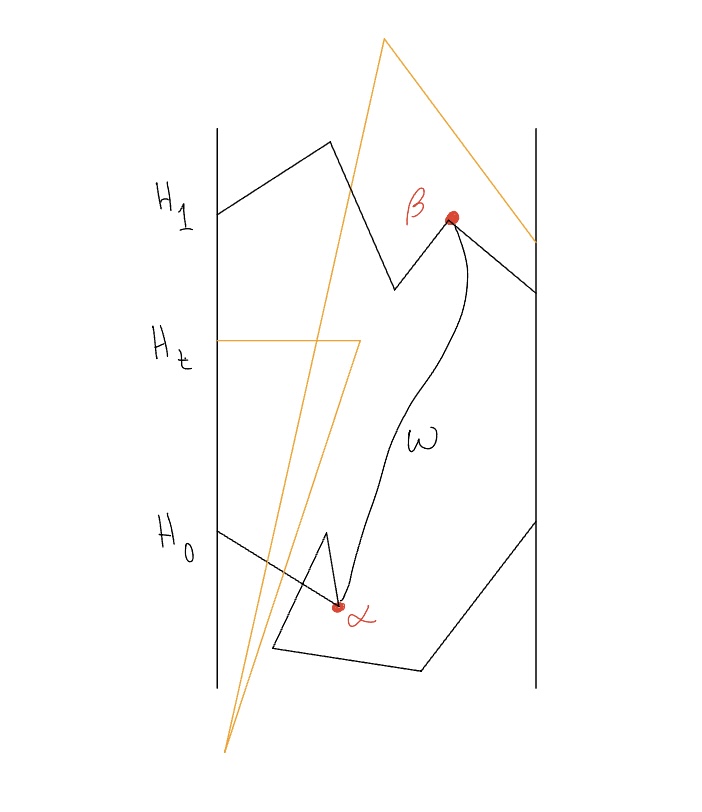}
\caption{A picture of the \textit{negation} of our simplifying assumption. In the general situation, it is possible that some of the pleated surfaces in the sweep-out might miss the image of $\omega$ completely.}
\label{fig:negation}
\end{figure}

\textit{Simplifying assumption:} Let us assume, without justification, that for each $t$, $\omega$ intersects the pleated surface $H|_{t}$-- see Figure \ref{fig:negation}.   

\vspace{2 mm}

\begin{figure}
\centering
\vspace{1cm}
\includegraphics[scale = 0.2]{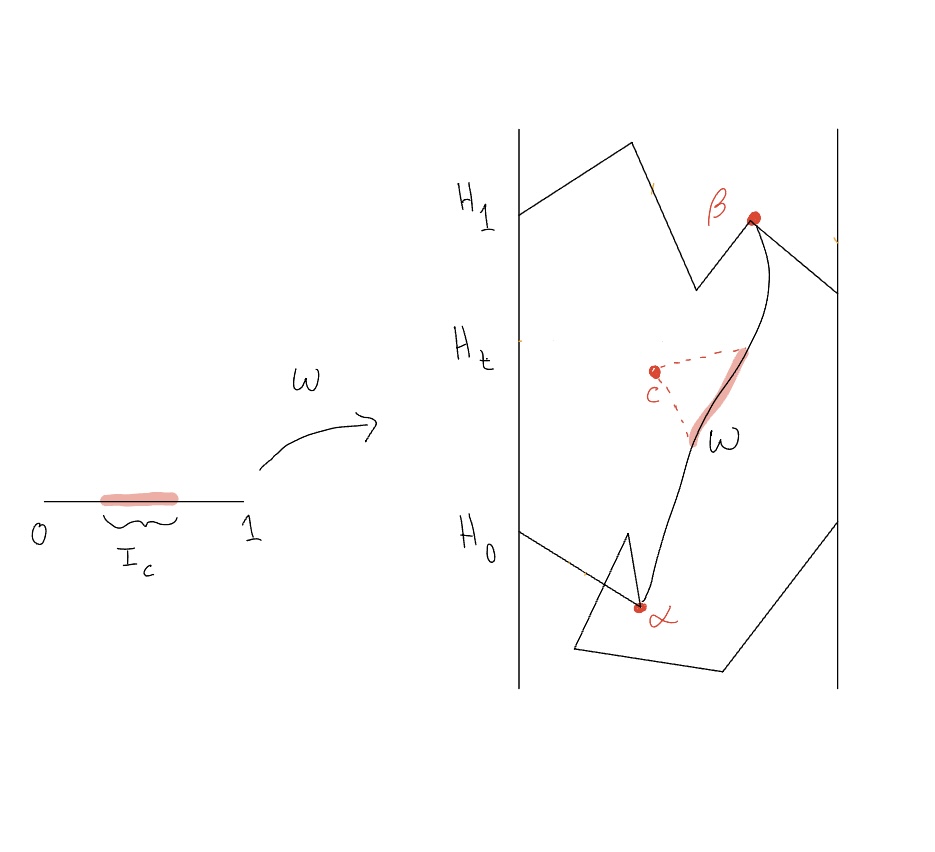}
\caption{A schematic cartoon of $I(c)$. The red dot denotes the location of the geodesic representative for $c$ in the ambient $3$-manifold.}
\label{fig:I(c)}
\end{figure}

For each curve $c$, let $I(c) \subset [0,1]$ denote the set of times $t$ for which $c$ admits a representative contained in $H|_{t}$ and based at a point in $\omega \cap H|_{t}$, with length at most $L$, for some value of $L$ chosen to guarantee that the various intervals $I(c)$ completely cover $[0,1]$-- see Figure \ref{fig:I(c)}. Note that this is possible only because of our simplifying assumption! 

In any case, the crucial point is that whenever $I(c) \cap I(c') \neq \emptyset$, then $c$ and $c'$ are \textit{simultaneously short} on some hyperbolic surface, and so a simple application of the collar lemma implies that their distance in $\mathcal{C}(S)$ is uniformly bounded\footnote{Each time $c$ intersects $c'$, it must cross through the its embedded tubular neighborhood and therefore one obtains a bound on $i(c,c')$ in terms of the lengths of $c$ and $c'$.}. It follows that (at least under the assumption that $\alpha$ has a sufficiently short geodesic representative in $M$), the sequence $\left\{x_{0}= \alpha, x_{1},..., x_{n}= \beta \right\}$-- where $I(x_{i}) \cap I(x_{i+1}) \neq \emptyset$ -- represents a path with uniformly bounded jumps in $\mathcal{C}(S)$ from $\alpha$ to $\phi(\alpha)$. 

So, a lower bound on $\ell_{\delta}(\omega)$ in terms of $d_{\mathcal{C}}(\alpha, \phi(\alpha))$ follows from a lower bound on $\ell_{\delta}(\omega)$ in terms of $n$. This follows readily from the standard fact that given $L$ and $\delta$, there is a uniform upper bound $B$ on the number of closed geodesics with length at most $L$ passing through a given point $p \in M_{\delta}$. Indeed, one imagines breaking $\omega$ up into subsegments of length at most $1$, and then arguing that each such subsegment meets at most some $B'$ curves of length $L+2$, and the key observation is that each $x_{i}$ is among one of these at most $B'$ curves for at least one of these subsegments. Thus, $n$ is at most $B'$ times the number of such subsegments, which is obviously proportional to $\ell_{\delta}(\omega)$, as required.

\subsubsection{Effective versions:} In joint work of the author with Patel and Taylor \cite{APT}, we obtain an \textit{effective} version of this theorem in the sense that we estimate the growth of $K$ as a function of $\chi(S)$ \cite{APT}. In particular, for a fixed $\delta$, we prove that $K$ grows polynomially in $|\chi(S)|$ of degree at most $30$. The most difficult feature of the argument is figuring out how to handle the possibility that some time-slices of the sweep-out don't intersect $\omega$, i.e. one needs to consider the possibility that Figure \ref{fig:negation} actually occurs and that therefore our simplifying assumption does not hold. 

\subsection{Systole length} \label{Systole finite type}

The \textit{systole}\index{systole} of a Riemannian manifold $M$ is the (or, perhaps ``a'') shortest closed geodesic. The \textit{systole length} of $M$, denoted $\textnormal{sys}(M)$, is the length of a systole. Any $M_{\phi}$ as described above will have at least one systole and therefore a well-defined and positive systole length. Minsky's bounded geometry theorem relates the systole length of $M_{\phi}$ to subsurface distances of the attracting and repelling laminations $\lambda^{+}, \lambda^{-}$ of $\phi$ \cite{Minsky}:

\begin{theorem} \label{thm:systole}  Let $\widetilde{M}$ be a hyperbolic $3$-manifold homeomorphic to $S \times \mathbb{R}$ (for example, $\widetilde{M}$ could be the infinite cyclic cover of some $M_{\phi}$ described above). Then given $\delta$, there is $J$ depending only on $\delta$ and $\chi(S)$, so that if $\alpha \in \mathcal{C}^{0}(S)$ has geodesic length at most $\delta$ in $M_{\phi}$, then 
\[ d_{Y}(\lambda^{-}, \lambda^{+}) > J, \]
for some subsurface $Y \subset S$ with $\alpha$ in its boundary. 
Conversely, if $d_{Y}(\lambda^{-}, \lambda^{+}) > J$, then $\ell(\partial Y) < \delta$. 
\end{theorem}

\subsubsection{Geometric intuition:} 
\begin{figure} \label{hierarchy}
\centering
\vspace{1cm}
\includegraphics[scale=0.2]{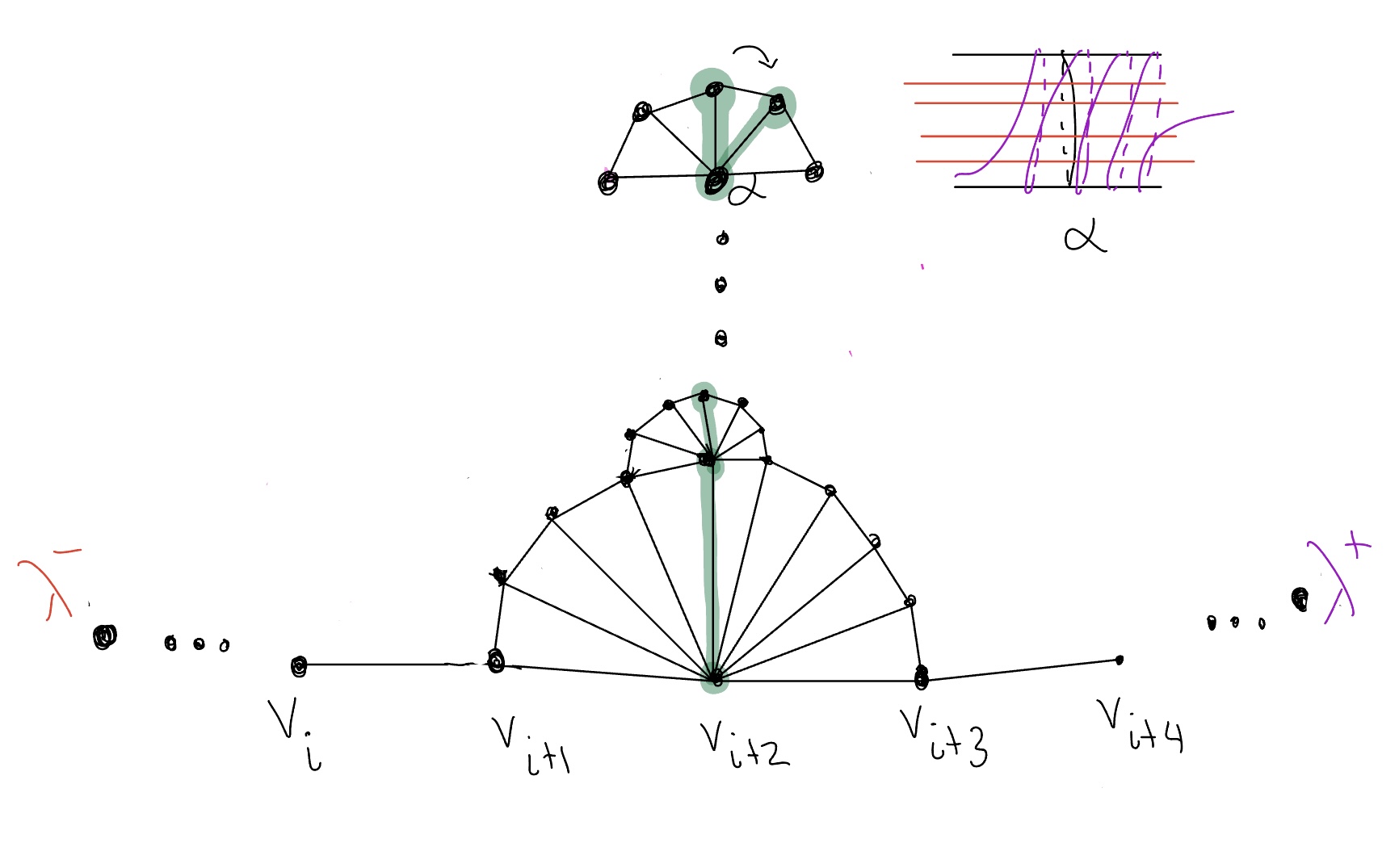}
\caption{A schematic drawing of a hierachy between a pair of ending laminations $\lambda^{-}, \lambda^{+}$. The intuitive idea is that a hierarchy is comprised first of a geodesic in the curve complex $\mathcal{C}(S)$ (pictured near the bottom), and then for each vertex $v_{i+2}$ along this geodesic, one connects the immediately preceding and immediately following vertices \textit{in the curve complex of the subsurface corresponding to the complement of $v_{i+2}$}. As one proceeds higher up in the hierarchy, the geodesics there live in curve complexes of subsurfaces that have smaller and smaller topological complexity, by virtue of being complementary to multi-curves with more and more components. At the very``top'' of the hierarchy are geodesics in annular complexes. Each geodesic has a length that is proportional to the subsurface projection distance between $\lambda^{-}$ and $\lambda^{+}$-- we illustrate this with a potential picture for how $\lambda^{\pm}$ might look within a given annulus. The distance between their projections in that particular annular complex is roughly $4$, and this is also the length of the corresponding geodesic in the hierarchy. The vertical green edges represents a so-called \textit{slice} of the hierarchy, which in this informal discussion the reader can think of as being a multi-curve comprised of vertices appearing somewhere in the hierarchy, together with a choice of a transverse curve to the core $\alpha$ of the annulus at the top of the slice. A move in a resolution is pictured, in which the slice ``clicks'' forward by moving one unit along the geodesic in the annular complex for $\alpha$, forming a new slice in which the top vertex is replaced with the one immediately to its right. As the resolution moves from slice to slice, the projection of that slice to a given subsurface resembles $\lambda^{+}$ more and more, and $\lambda^{-}$ less and less.}
\label{fig:hierarchy}
\end{figure}

To get a grasp for how subsurface projections relate to geodesic length, we consider the hierarchy and model-manifold machinery of Masur--Minsky and Brock--Canary--Minsky (\cite{MasurMinskyII}, \cite{BCM})\footnote{We remark that the description we are about to give is ahistorical, since Minsky's bounded geometry theorem was proved before the model manifold machinery was formalized, and is in fact an \textit{ingredient} for doing so.}. Without getting into too many of the formal details, the crucial point is that given a hyperbolic $3$-manifold $M$ homeomorphic to $S \times \mathbb{R}$ and with ending laminations $\lambda^{-}, \lambda^{+}$  there is a combinatorial object called a \textit{hierarchy}\index{hierarchy} built out of paths in the curve complexes of various subsurfaces of $S$ with the following properties (see the poorly drawn but hopefully easy to understand Figure \ref{fig:hierarchy}): 
\begin{enumerate}
\item The length of the path in $\mathcal{C}(Y)$ is roughly equal to $d_{Y}(\lambda^{-}, \lambda^{+})$.
\item There is a natural way (or really, a collection of ways)-- called a \textit{resolution}-- to order the vertices occurring in the geodesic paths of the hierarchy, and then to navigate through the various paths one after the other. 
\item For each $Y \subseteq S$, Brock-Canary-Minsky associate a so-called \textit{block} which is a metric manifold homeomorphic to $Y \times [0,1]$, and they use the data of the hierarchy to give a way of gluing the blocks together to form a manifold homeomorphic to $S \times \mathbb{R}$ which is bi-Lipschitz equivalent to $M$. Since the geometry of this manifold models that of $M$, it is referred to as the \textit{model manifold} for $M$; moving efficiently from one end of the model to the other corresponds in a precise way to moving through a resolution of the hierarchy.  
\end{enumerate}

Now, if one is moving through $M$ from $\lambda^{-}$ to $\lambda^{+}$, it would be very inefficient to enter deeply into a Margulis tube associated to a very short geodesic $\alpha$. Instead, one would move around the boundary of the tube, which means that, in the model manifold, one is traveling within a block associated to $S \setminus \alpha$. Intuitively, the shorter $\alpha$ is, the longer the required detour will be through the block for $S \setminus \alpha$ since in the limiting case, the length of $\alpha$ is $0$ and so $\alpha$ corresponds to a cusp and then the block for $S \setminus \alpha$ is actually the entire manifold. Thus, using property $1$ and $3$ above, the subsurface projection $d_{S \setminus \alpha}(\lambda^{-}, \lambda^{+})$ should be inversely proportional to $\ell_{M}(\alpha)$.

Recently, Viaggi obtained an effective version of the converse direction of Theorem \ref{thm:systole} \cite{Vi}, which essentially states the existence of some universal and computable constant $b$ and polynomials $p_{1}, p_{2}$ of degrees $345$ and $348$, respectively, so that if the projection $d_{Y}(\lambda^{-}, \lambda^{+})$ is at least $p_{1}(|\chi(Y)|)$, then $\ell(\partial Y)$ is at most 
\[ \frac{p_{2}(|\chi(Y)|)}{d_{Y}(\lambda^{-}, \lambda^{+}) - b \log|\chi(S)|}. \]

\subsection{Volume} \label{Volume finite type}

The hyperbolic volume of $M_{\phi}$ is coarsely captured by the action of $\phi$ on the pants graph\index{pants graph}, as described in the following theorem due to Brock \cite{Brock2}:

\begin{theorem} \label{thm:brock1}
There is a constant $L$ depending only on $\chi(S)$ so that 
\[ \tau_{\mathcal{P}}(\phi) \sim_{L} \textnormal{Vol}(M_{\phi}). \]
where $\tau_{\mathcal{P}}(\phi)$ denotes the asymptotic translation length of $\phi$ in the pants graph. 
\end{theorem}

\subsubsection{Geometric intuition:} Given a path $\left\{p_{0},..., p_{n} = \phi(p_{0}) \right\}$ in the pants graph, there is a corresponding sequence of curves $\alpha_{1},..., \alpha_{n}$, where $\alpha_{i}$ is the curve in $p_{i}$ that replaces a curve in $p_{i-1}$ through an elementary move. We can then consider the hyperbolic $3$-manifold $\widehat{M_{\phi}}$ obtained from $M_{\phi}$ by drilling out each $\alpha_{i}$. In \cite{Agol}, Agol describes how to view $\widehat{M_{\phi}}$ as being built from $n$ regular ideal octahedra, each corresponding to one of the elementary moves in the provided path. An upper bound on the volume of $\widehat{M_{\phi}}$ in terms of $n$ follows readily, and then one obtains an upper bound on the volume of $M_{\phi}$ itself with a theorem of Thurston that states that the volume of a drilled manifold is always bounded below by the original (see either \cite{CM} or chapter 6 of \cite{Thurston4}).

As is generally the case, lower bounds on geometric quantities are harder and more complicated to obtain. But again, the key tool is sweep-outs: one sweeps through $M$ with pleated surfaces and keeps track of short curves along the way. First, one bounds from below the number of relatively short curves by the pants graph distance (using the hierarchy machinery of Masur-Minsky), and then one uses the standard fact that each short curve in a hyperbolic $3$-manifold contributes definitely to its total volume.

\subsubsection{Effective versions:} The argument for the upper bound on volume sketched directly above is already effective. As for the lower bound, work of the author with Taylor and Webb establishes a bound on $L$ that depends roughly factorially on $\chi(S)$ \cite{ATW}.

\subsection{Weil-Petersson translation length} \label{WP}

A second result of Brock coarsely relates the $\tau_{\mathcal{P}}(\phi)$-- and by the previous subsection, also the volume of $M_{\phi}$-- to the Weil-Petersson translation length $\tau_{WP}(\phi)$ of $\phi$ acting on the Teichm{\"u}ller space \cite{Brock1}:

\begin{theorem} \label{thm:WP}
There is some $L' = L'(S)$ so that 
\[ \tau_{\mathcal{P}}(\phi) \sim_{L'} \tau_{WP}(\phi). \]

\end{theorem}

\subsubsection{Geometric Intuition:} Given $X,Y \in \mathcal{T}(S)$, choose pants decompositions $P_{X}, P_{Y}$ that are shortest possible on $X,Y$, respectively. The idea behind Theorem \ref{thm:WP} is to prove that the map 
\[ \psi: \mathcal{T}(S) \rightarrow \mathcal{P}^{0}(S), \psi(X) := P_{X}\]
is a $\mbox{Mod}(S)$-equivariant quasi-isometry. The statement about $\tau_{WP}(\phi)$ then follows readily as a corollary.

To see that $\psi$ is (coarsely) Lipschitz, consider a pair of metrics $X,Y$ that are distance $1$ in the Weil-Petersson metric.

Another fundamental result in Weil-Petersson geometry is the fact that the Weil-Petersson distance is bounded above by Teichm{\"u}ller distance \cite{Linch}, and therefore $d_\mathcal{T}(X,Y) \leq 1$ as well. We can now apply a well-known inequality of Wolpert \cite{Wolpert2}: for any closed curve $\alpha$ on $X$, one has that 
\[ \ell_{X}(\alpha) \leq \exp(d_{\mathcal{T}}(X,Y)) \cdot \ell_{Y}(\alpha). \]
From this, and from the fact that $P_{Y}$ has length at most the Bers constant $B(S)$ on $Y$, it is straightforward to see that $P_{Y}$ has length at most some constant $B'= B'(S)$ on $X$. The collar lemma now implies a uniform upper bound on the quantity $i(P_{X}, P_{Y})$; indeed, since both have length at most $B'$ on the same surface, each component has a definitely long embedded collar and so $P_{Y}$ can cross $P_{X}$ only so many times before growing very long. 

Finally, one simply needs to argue that $d_{P}(P_{X}, P_{Y})$ can be bounded above in terms of the geometric intersection number $i(P_{X}, P_{Y})$. This follows immediately from the fact that there are only finitely many $\mbox{Mod}(S)$-orbits of pairs of pants decompositions with a fixed intersection number. 

The reverse bound--namely, that $\psi$ does not contract distances an arbitrary amount--is more complicated and involves some of the subsurface projection and hierarchy machinery of Masur-Minsky mentioned above. But at the very least, it is easy to deal with one potential source for large contractions: the set of all points in $\mathcal{T}(S)$ for which the same pants decomposition $P$ is shortest will all map to $P$ by $\psi$. A key geometric property of the Weil--Petersson metric is that it is incomplete: the Weil--Petersson distance $X$ to the degenerate hyperbolic surface obtained by pinching each component of $P_{X}$ to points, is uniformly bounded and depends only on the topology of $S$ (see for instance \cite{Wolpert1}). 

\subsubsection{Effective versions.} The work of the author joint with Taylor and Webb \cite{ATW} also establishes an effective relationship between Weil-Petersson translation length and pants distance. 

A consequence of this and the last subsection is that the volume of $M_{\phi}$ is coarsely equal to the Weil--Petersson translation length, and one could ask whether there are effective versions of this relationship that do not factor through the pants graph. Indeed, Kojima--McShane \cite{KJ} and Brock-Bromberg \cite{BB2} are able to obtain very concrete and effective upper bounds on the volume in terms of Weil--Petersson translation length, as follows: 
\[ \mbox{Vol}(M_{\phi}) \leq \frac{3}{2} \sqrt{2 \pi |\chi|} \cdot \tau_{WP}(\phi).\]

\subsection{Cusp area}

Assuming that $S$ has at least one puncture (and no boundary for simplicity), one can consider the action of $\phi$ on the arc complex $\mathcal{A}(S)$\index{arc complex}. Moreover, the puncture will give rise to a rank-$2$ cusp in $M_{\phi}$, and by a \textit{maximal} neighborhood $C$ of this cusp we will mean the largest embedded tubular neighborhood of it. The boundary $\partial C$ of a maximal neighborhood is naturally isometric to a Euclidean torus. A theorem of Futer and Schleimer relates the area of this torus to the $\mathcal{A}(S)$ translation length of $\phi$ \cite{FS}:

\begin{theorem} \label{thm:FS} 
When $S$ is singly punctured, 
\[ \frac{\tau_{\mathcal{A}}(\phi)}{450 \chi(S)^{4}} \leq \textnormal{Area}(\partial C) \leq 9 \chi(S)^{2} \tau_{\mathcal{A}}(\phi),\]
where $\tau_{\mathcal{A}}(\phi)$ denotes the asymptotic translation length of $\phi$ on $\mathcal{A}(S)$. 
\end{theorem}

\subsubsection{Geometric intuition:} To understand the upper bound on cusp area, begin with a path $\left\{a_{0},..., a_{n} = \phi(a_{0}) \right\}$ in the arc complex. Since $a_{i}$ and $a_{i+1}$ are disjoint, they live together within some ideal triangulation $T_{i}$ of $S$. We can then consider a pleated surface $X_{i}$ in $M_{\phi}$ with pleating locus equal to $T_{i}$. Futer and Schleimer prove Theorem \ref{thm:FS} by analyzing the intersection of $\partial C$ with $X_{i}$, and by arguing that it must correspond with a polygonal curve on $X_{i}$ with length bounded above by an explicit linear function of $\chi(S)$. They then estimate the distance between $a_{i} \cap \partial C$ and $a_{i+1} \cap \partial C$ and show that it is also at most some explicit linear function of $\chi(S)$. One then imagines $\partial C$ as being covered by (potentially immersed) cylindrical pieces with circumference at most the length of the aforementioned polygonal curve and height the distance between $a_{i} \cap C$ and $a_{i+1} \cap C$, and the desired bound follows. Unsurprisingly, lower bounds on these geometric quantities are again obtained through sweep-out arguments. 

Note that Theorem \ref{thm:FS} is already effective, although it is unclear how sharp it is. 

\section{Infinite type surfaces} \label{infinite type}

In the last decade, there has been a flurry of new work aimed at generalizing the geometric tools described above for studying the mapping class group of a finite type surface, to those (orientable) surfaces with fundamental groups that are not finitely generated. 

We refer the reader to Aramayona--Vlamis \cite{AV} for a detailed exposition of the theory of infinite type surfaces. For our purposes, it suffices to note that up to homeomorphism, such manifolds are determined by the same three invariants as in the finite type case ($g,b$, and $p$), together with a fourth invariant-- the homeomorphism type of the space of ends. 

Perhaps unsurprisingly, much of the technology for studying the mapping class group of a finite type surface does not carry over to the infinite type setting. For example: 

\begin{itemize}
\item For a finite type surface $S$, $\mbox{Mod}(S)$ is generated by finitely many Dehn twists. For many infinite type surfaces, it is not even true that $\mbox{Mod}$ is \textit{topologically} generated by Dehn twists (the closure of the set of all mapping classes generated by twists is proper)-- see for instance Patel-Vlamis \cite{PV}. 
\item As a consequence of residual finiteness and finite generation, the mapping class group of a finite type surface is Hopfian\index{hopfian} (every onto homomorphism from the group to itself is an automorphism); this is false, even for fairly simple infinite type surfaces \cite{AV}. 
\item In the finite type case, the Nielsen realization theorem states that finite subgroups of $\mbox{Mod}(S)$ are precisely the groups that can be realized as isometry groups of hyperbolic metrics on $S$. In infinite type, only one direction is true: finite subgroups of the mapping class group are realizable as isometry groups (see for instance \cite{ACCL}), but on the other hand, there are infinite type surfaces for which \textit{any countable group} can be realized as the isometry group for some hyperbolic metric \cite{APV}. 
\item $\mbox{Mod}(S)$ is acylindrically hyperbolic\index{acylindrically hyperbolic} if and only if $S$ has finite type \cite{BG}. 
\end{itemize}
The list goes on. For example, as we saw in the previous sections, the curve and arc complexes are crucial tools for probing the structure of the mapping class group, and almost always the starting point for using either of them towards this end is the fact that they are infinite diameter and Gromov hyperbolic. On the other hand, it is straightforward to see that if one defines $\mathcal{C}(S)$ and $\mathcal{A}(S)$ in the same way for an infinite type surface, their diameters are $2$, and this short-circuits any hope of porting over some of the most interesting applications to this more complicated setting. 

For these reasons, it is all the more interesting when finite type results \textit{do} carry over. Even more fascinating are those results which are stated for infinite type surfaces, but which imply highly non-trivial results about \textit{finite} volume hyperbolic $3$-manifolds. The connection here comes from taut foliations\index{taut foliation}; given a taut, depth-one foliation of a closed hyperbolic $3$-manifold, a non-compact leaf will be an infinite type surface with finitely many ends each accumulated by genus. See for instance Calegari's book \cite{C} on foliations for more details. 

Such a leaf is dense in an open submanifold that is naturally a mapping torus for a homeomorphism on an infinite type surface. And this homeomorphism will always be of a certain, well-behaved type that we explore in more depth in the next subsection.

\subsection{Endperiodic maps} 
The idea of an endperiodic homeomorphism $f: S \rightarrow S$ is that it is a map which, at least up to taking (potentially negative) powers, sends each end into itself. We refer the reader to Section 2.2 of \cite{FKLL} for a more detailed description of the theory we only touch on below; see also Cantwell-Conlon-Fenley \cite{CCF}. 

\begin{definition} \label{def:endperiodic} A homeomorphism $f: S \rightarrow S$ is \textit{endperiodic}\index{endperiodic} if there is some natural number $n \in \mathbb{N}$ so that for each end $E$ of $S$ admits a neighborhood $U_{E}$ so that either $f^{n}(U_{E}) \subsetneq U_{E}$ or $f^{-n}(U_{E}) \subsetneq U_{E}$, and so that $\left\{f^{nm}(U_{E}): m \in \mathbb{N} \right\}$ (resp. $\left\{f^{-nm}(U_{E}) : m \in \mathbb{N} \right\}$) forms a neighborhood basis of $E$.  
\end{definition}

An end $E$ for which $f^{n}(U_{E}) \subsetneq U_{E}$ is called \textit{attracting}; when $f^{-n}(U_{E}) \subsetneq U_{E}$, $E$ is said to be \textit{repelling}. If for each attracting (resp. repelling) $U_{E}$, one has $f(U_{E}) \subset U_{E}$ (resp. $f^{-1}(U_{E}) \subset U_{E}$), then the set of neighborhoods $\left\{U_{E} \right\}$ is said to be a collection of \textbf{tight nesting neighborhoods}.

\begin{definition} An endperiodic homeomorphism $f$ is \textbf{irreducible} if no curve is sent to itself by any non-trivial power of $f$; for any line $\ell$ (a topological embedding of $\mathbb{R}$) in $S$ with one endpoint in an attracting end and the other in a repelling end, $\ell$ is never sent to itself by any power of $f$; and there is no curve $\gamma$ and integers $n,m \in \mathbb{Z}$ ($n> m$) so that $f^{n}(\gamma)$ is contained in some $U_{E}$ as in Definition \ref{def:endperiodic} for $E$ an attracting end, and $f^{m}(\gamma)$ is contained in some $U_{E'}$ for $E'$ a repelling end. 
\end{definition}

The definition of irreducibility is designed to trigger the celebrated hyperbolization results of Thurston for the mapping torus $M_{f}$; for example, the absence of curves that are sent to themselves by powers of $f$ help to rule out the existence of $\pi_{1}$-injective tori. Indeed, Proposition $3.1$ of Kim--Field--Leininger--Loving states that when $f$ is endperiodic and irreducible on a surface $S$ with finitely many ends each accumulated by genus, then $M_{f}$ is the interior of a compact, irreducible, atoroidal $3$-manifold $\overline{M_{f}}$ with incompressible boundary \cite{FKLL}. It follows that $\overline{M_{f}}$ admits a convex hyperbolic metric. 

\subsection{Volume.} 

Putting together Field--Kim--Leininger-Loving and Field--Kent--Leininger--Loving (\cite{FKLL}, \cite{FKeLL}), one obtains a complete analog of Brock's theorem for the finite-type setting that relates pants distance to volume of the mapping torus: 

\begin{theorem} \label{thm:volume infinite} Let $f: S \rightarrow S$ be an irreducible and endperiodic homeomorphism on a surface $S$ with finitely many ends, each accumulated by genus. Let $\mbox{Vol}(\overline{M_{f}})$ denote the infemal volume for all convex hyperbolic metrics on the compact $3$-manifold $\overline{M_{f}}$. Then there are constants $L(f)$ and $K$ so that 
\[L(f) \cdot \tau_{\mathcal{P}}(f) \leq \mbox{Vol}(\overline{M_{f}}) \leq K \cdot \tau_{\mathcal{P}}(f), \]
where: 
\begin{itemize}
\item $\tau_{\mathcal{P}}(f)$ is the asymptotic translation length of $f$ on the subgraph of $\mathcal{P}(S)$ spanned by pants decompositions $P$ so that for some finite $n$, $P$ and $f^{n}(P)$ differ by only finitely many elementary moves. 
\item $L(f)$ is a constant depending only on the \textbf{capacity} $C(f)$ of the homeomorphism $f$, which measures the topological complexity of both $\partial \overline{M_{f}}$, and the minimal finite type subsurface of $S$ that acts as a \textbf{core} for the map $f$ (a compact subsurface $Y$ for which $S \setminus Y$ is a collection of tight nesting neighborhoods); and 
\item $K$ is a universal constant which can be taken to be equal to the volume of a regular right-angled ideal octahedron.
\end{itemize}

\end{theorem}

We note the differences and similarities between Theorems \ref{thm:volume infinite} and \ref{thm:brock1}; first, the upper bound on volume is essentially a perfect parallel to the finite type setting. On the other hand, the lower bound changes from being stated in terms of a constant depending only on the topology of the fiber, to a constant that now depends on the dynamics of the \textit{map} $f$. Upon brief reflection, this is not so surprising: the fiber now has infinitely generated fundamental group, so one would expect that any \textit{direct} analog of Brock's theorem-- proven in the way that Brock proves it-- would involve an infinitely large multiplicative error term. One of the elegant work-arounds used to prove Theorem \ref{thm:volume infinite} is to focus in on a pair of finite type surfaces -- the core of $f$ and the boundary of $\overline{M_{f}}$-- that, while being very small in comparison to $S$ itself, capture enough of the dynamics of $f$ to encode geometric properties of $\overline{M_{f}}$.  

We also point out that while there are many genuinely new and difficult ideas introduced by the authors to prove this theorem, the essential strategy is similar to that of Agol \cite{Agol} and Brock (\cite{Brock1}, \cite{Brock2}): for the upper bound, one relates each elementary move in the pants graph to a ``block'' sub-manifold of definite volume in a drilled out version of $\overline{M_{f}}$. For the lower bound, one imagines a sweep-out by pleated surfaces, tries to relate pants distance to the number of short curves that occur in the sweep-out, and then appeals to the fact that each short curve contributes definitely to volume. To establish that either of these basic outlines actually works requires extreme care and ingenuity from the authors, especially the lower bound because a priori one loses all of the usual control over the hyperbolic geometry when the cross sections of a sweep-out are themselves infinite type surfaces.

\subsection{Systoles} 

We next discuss an infinite type analog of Minsky's theorem \ref{thm:systole} that inversely relates lengths of closed geodesics in a hyperbolic manifold homeomorphic to $S \times \mathbb{R}$ to the distance in the curve complex of its complement, between the projections of the ending laminations. For the same reasons discussed in the previous subsection, one would expect that obtaining such a theorem for pseudo-Anosov-like homeomorphisms in an infinite type setting by following the same rough outline of the original proof, would be extremely challenging. For starters, the complement of a curve is itself an infinite type subsurface whose curve complex is therefore of finite diameter. So care needs to be taken even in the way that such a theorem should be stated; for example, instead of a statement that applies to curves, perhaps one could formulate a theorem that pertains to a given compact subsurface and which relates the length of its boundary to the projections to its (infinite diameter) curve complex of the appropriate laminations.  

Morover, one needs a theory of ending laminations for the homeomorphisms on an infinite type surface $S$ which we imagine as being the analogs of pseudo-Anosovs. A natural candidate would be the so-called \textit{Handel--Miller} laminations associated with any endperiodic map $f:S \rightarrow S$ on an infinite type surface with finitely many ends, each accumulated by genus-- we refer to \cite{CCF} for the details of this construction. 

In any case, this is indeed the approach taken by Whitfield \cite{Whitfield}, who proved the following:

\begin{theorem} \label{thm:Brandis} For any $D, \epsilon >0$, there is $K$ so that if $f:S \rightarrow S$ is irreducible and endperiodic with the property that the boundary of $\overline{M_{f}}$ consists solely of genus $2$ surfaces and so that the capacity $C(f) \leq D$, then for any compact essential subsurface $Y \subset S$, one has 
\[ d_{Y}(\Lambda^{-}, \Lambda^{+}) \geq K \Rightarrow \ell(\partial Y) \leq \epsilon, \]
where:
\begin{itemize}
\item $\Lambda^{\pm}$ are the Handel--Miller laminations for $f$ on $S$, and 
\item $\ell$ is the length of $\partial Y$ in the unique hyperbolic metric on $\overline{M_{f}}$ with totally geodesic boundary. 
\end{itemize}

\end{theorem}

The core of the idea behind Whitfield's strategy is to apply a theorem of Landry-Minsky-Taylor \cite{LMT} to double $\overline{M_{f}}$ over its boundary components to get a closed hyperbolic $3$-manifold $N$, and then to apply Minsky's original theorem to this object. By keeping careful track of the dynamics of the suspension flow induced by $f$ on $M_{f}$, one can do this in such a way that flow lines glue up smoothly to produce a $1$-dimensional foliation on $N$. One then observes that if there is a closed surface $\Sigma$ in the double transverse to this flow and meeting each flow line once, then $N$ is necessarily fibered over the circle, with fiber isotopic to $\Sigma$. 

In the context of Whitfield's argument, she carefully selects a core for the homeomorphism $f$ on $S$; its double will be closed and will satisfy the requirements outlined in the previous paragraph, such that it arises as a fiber of a fibration on the double. Therefore, $N$ is fibered by a surface whose genus is bounded in terms of the topological complexity of the core chosen on $S$. There must also be some pseudo-Anosov monodromy on $N$, with ending laminations $\lambda^{+}, \lambda^{-}$.  

Now, if $Y \subset S$ is some compact subsurface with the property that $d_{Y}(\Lambda^{-}, \Lambda^{+}) \ge 3$, one can argue that $Y$ admits an embedding into the doubled manifold whose image is isotopic to some $W$ that lives in the double of the core of $S$. The devil is of course in the details, but the intuition is that the core captures all of the interesting dynamics of $f$, and so a surface living outside of it should not witness a large projection. Furthermore, since the flow on the double was carefully constructed from the flow on $M_{f}$, Whitfield also shows that $d_{Y}(\Lambda^{-}, \Lambda^{+})$ can be uniformly bounded above by $d_{W}(\lambda^{-}, \lambda^{+})+ 2$. Therefore,
\[ d_{Y}(\Lambda^{-}, \Lambda^{+}) \ge K \Rightarrow d_{W}(\lambda^{-}, \lambda^{+}) \ge K-2.  \]

The pieces are now all in place: Minsky's theorem \ref{thm:systole} implies that the boundary of $W$ has short geodesic length in the hyperbolic metric on $N$. The only remaining detail is to relate this to the geodesic length in $\overline{M_{f}}$ of $\partial Y$, and this is where the assumption that each boundary component of $\overline{M_{f}}$ is a genus $2$ surface comes into play. Indeed, Whitfield chooses the doubling map to correspond to the hyperelliptic involution on each totally geodesic genus $2$ boundary surface, and then argues that with this choice of gluing, the unique hyperbolic structure with totally geodesic boundary of the compactified mapping torus isometrically embeds in the double. 

As alluded to at the beginning of Section \ref{infinite type}, in the opinion of the author, the most interesting theorems in the infinite type realm are those which are not only inspired by finite type theorems, but which produce new corollaries in the finite type setting. Whitfield's result achieves this standard \cite{Whitfield}: 

\begin{corollary} \label{cor:Brandis} For any $\epsilon \ge 0, g \ge 2$, there is a closed and fibered hyperbolic $3$-manifold $N$ which admits a totally geodesic surface $\Sigma$ of genus $g$ whose systole has length at most $\epsilon$. 
\end{corollary}

\section{Fixed points} \label{Fixed points}

We end this survey with a description of recent work due to the author, joint with Futer and Taylor \cite{AFT}, and which pertains to an invariant of a pseudo-Anosov that we have yet to discuss: the number of fixed points. Of course, a pseudo-Anosov corresponds to an entire isotopy class of homeomorphisms, but one can prove that the minimum number of fixed points, taken over the entire pseudo-Anosov homotopy class, is realized by the representative $f'$ described in Subsection \ref{subsec:pA}. 

Intuitively, the number of fixed points should correlate with the complexity of the map. There are several ways to explain this heuristic-- we will do so using train tracks: for each pseudo-Anosov $f$, there is an invariant train track\index{train track} $\tau$ on the surface carrying its attracting lamination. As one iterates $f$, the orbits of $\tau$ are carried by $\tau$ itself, in such a way where $f^{n}(\tau)$ wraps potentially many times around $\tau$. Each time an edge of $\tau$ is mapped over itself, a simple application of the Brower fixed point theorem shows that there must be a corresponding fixed point. The dynamics of how the edges of $\tau$ are sent to train-paths on $\tau$ is encoded by an irreducible matrix $A_{f,\tau}$ whose $(i,j)^{th}$ entry corresponds to the number of times the image of the $i^{th}$ branch of $\tau$ under $f$ passes over the $j^{th}$ edge, and so the fixed points of $f^{n}$ grow like the trace of $A^{n}$.  

Therefore, letting $\left\{\mbox{Fix}(f^{n}) \right\}$ denote the sequence of the number of fixed points of $f^{n}$ as we take higher and higher powers, standard Perron-Frobenius theory dictates that it should grow exponentially with base given by the highest eigenvalue of $A_{f,\tau}$, which is none other than $\lambda$, the dilatation of $f$. In turn, $\lambda$ is precisely the logarithm of the translation length $\ell_{T}(f)$ of $f$ on the Teichm{\"u}ller space with respect to the Teichm{\"u}ller metric:
\begin{equation} \label{eq:powers}
 \mbox{Fix}(f^{n}) \asymp \lambda^{n} \Rightarrow \log(\mbox{Fix}(f^{n})) \asymp n \cdot \ell_{T}(f). 
\end{equation}

Since there is a uniform lower bound on the dilatation of a pseudo-Anosov given solely in terms of the Euler characteristic of the underlying surface, one can deduce uniform lower bounds on the number of fixed points for a pseudo-Anosov $f^{n}$ on $S$ in terms only of $n$, and furthermore, this number grows exponentially in $n$. It is then natural to ask for what one can say more generally about maps that are not proper powers of any other map, but which are still ``dynamically complex'', for example: perhaps $f$ has a large translation length in one of the combinatorial complexes highlighted above. 

Such an assumption is insufficient, at least for certain complexes. Indeed, translation length can be very large in the pants complex without guaranteeing the existence of many fixed points. Indeed, as the example below demonstrates, there are fixed point-free pseudo-Anosovs that send an isotopy class to one that is disjoint from it, and this can occur simultaneously with large pants graph translation length.

One could instead focus on the curve complex and try to relate $\tau_{\mathcal{C}}(f)$ to $\mbox{Fix}(f)$, but this also misses a part of the picture because there are pseudo-Anosovs with many fixed points but with small curve complex translation length -- we leave it as an exercise to the reader to construct examples along these lines.  

The optimal statement should therefore involve \textit{some} definite curve complex translation as an assumption (in order to avoid maps as described in the example below, but ultimately the lower bound on the number of fixed points in such a statement should depend on a quantity that can be arbitrarily large relative to curve complex translation. Equation \ref{eq:powers} suggests that \textit{Teichm{\"u}ller} translation length might be a good candidate-- this motivates the following theorem due to the author, joint with Futer and Taylor \cite{AFT}: 

\begin{theorem} \label{thm:AFT} Let $f$ be a pseudo-Anosov on a surface $S$ with $\chi(S)<0$ so that for all simple closed curves $\alpha$, 
\[ i(\alpha, f(\alpha))>0. \]
(Note that this is equivalent to requiring that the curve complex translation length is at least 2.) Then
\[ \log (\textnormal{Fix}(f)) \sim \ell_{T}(f).\]

\end{theorem}

We remark that the assumption that every simple closed curve intersects its image can be removed, at the cost of a more complicated conclusion that no longer depends on $\ell_{T}$. Instead, the bound will depend on the size of the subsurface projections $d_{Y}(\lambda^{-}, \lambda^{+})$ to all subsurfaces $Y$ that are ``displaced'' by $f$ in the sense that $\partial Y$ and $f(\partial Y)$ intersect essentially: 

\begin{theorem} \label{thm:AFT2} Let $f$ be a pseudo-Anosov on a surface $S$ with $\chi(S)<0$, and let $\mathcal{D}_{f}$ denote the set of $f$-orbits of subsurfaces $Y$ for which $Y$ and $f(Y)$ overlap. Then 
\[ \log(\textnormal{Fix}(f)) \sim \ell_{\mathcal{C}(S)}(f) + \sum_{[Y] \in \mathcal{D}_{f}} [d_{Y}(\lambda^{-}, \lambda^{+})] + \sum_{[A] \in \mathcal{D}_{f}} [\log(d_{A}(\lambda^{-}, \lambda^{+})], \]
where the first sum takes place over non-annular $f$-orbits and the second sum corresponds to $f$-orbits of annular subsurfaces. Moreover, in the above sum we take $\log(0)$ to be $0$. 

\end{theorem}

\subsection{Key ideas} To get a feel for the strategy and the core arguments that go into these theorems, it can be helpful to first appreciate an example that demonstrates the necessity of the hypothesis that every simple closed curve intersects its image: 

\textbf{Example} If $M$ is a hyperbolic $3$-manifold fibering over the circle with first Betti number at least $2$, then $M$ will fiber in infinitely many ways. These fibrations are organized by lattice points in $H^{2}(M)$, where a given fibration is represented by the homology class of the fiber. The \textit{Thurston norm}\index{Thurston norm} assigns to each homology class the minimum absolute value-- taken over all embedded representatives-- of the Euler characteristic. With respect to this norm, the unit ball is polyhedral with rational vertices. The cone over each face of this polyhedron is called a \textit{fibered face}, and each fibration appearing within it is transverse to a single pseudo-Anosov flow in $M$. As one appraoches the boundary of a fibered face, the Thurston norm blows up, as does the dilatation of the corresponding monodromy map. 

With all of this in mind, consider a pair of fibers $S_{1}, S_{2}$ in the same fibered face. Each is associated to a corresponding monodromy pseudo-Anosov $f_{1}, f_{2}$. The common pseudo-Anosov flow to which both $S_{i}$ are transverse has the property that $f_{i}$ represents the first return map to the fiber $S_{i}$. Moreover, a fixed point of $f_{i}$ corresponds to a periodic orbit of this flow intersecting $S_{i}$ exactly one time. 

One can then form a new fiber $S_{3}$ obtained by surgically resolving any intersections between $S_{1}, S_{2}$ while maintaining transversality to the flow; this fiber will be an embedded representative of the homology class $[S_{1}]+ [S_{2}]$, and will therefore also be located in the same fibered face. It follows that any periodic orbit of the flow intersects $S_{3}$ at least twice, and therefore its monodromy map, $f_{3}$, will be fixed point-free. 

Now, let $\gamma \subset S_{1} \cap S_{2}$ be a simple closed curve that is essential on both $S_{1}$ and $S_{2}$. When the intersections between these fibers are resolved, the resulting surface $S_{3}$ will now contain two homotopic copies of $\gamma$, and the first return map $f_{3}$ sends one to the other\footnote{We point out a subtlety here: these two curves do \textit{not} represent the same homotopy class on the fiber $S_{3}$, for if they did, $f_{3}$ would not be pseudo-Anosov.} Therefore, we see that $f_{3}$ does \textit{not} satisfy the key hypothesis of Theorem \ref{thm:AFT}.

So we see that some version of the hypothesis is necessary. With it, the strategy for proving Theorem \ref{thm:AFT} proceeds as follows: 

\begin{enumerate}
\item One first argues that the number of fixed points of $f$ is captured by the dynamics of the pseudo-Anosov representative: it relates to the number of times a so-called \textit{singularity-free saddle connection}\index{saddle connection} intersects its image under $f$. Concretely, a \textit{singularity-free rectangle} in the universal cover of $S$ equipped with the quadratic differential metric associated to $f$ is a rectangle with vertical and horizontal sides containing no singularities in its interior. A \textit{saddle connection} is a straight line segment connecting two singularities, and a \textit{singularity-free saddle connection} is a saddle connection spanning a singularity-free rectangle, in the sense that it connects two singularities on the boundary of that rectangle.

By a careful application of the contraction mapping theorem, one can relate the number of fixed points contained in a given singularity-free rectangle to the number of times a corresponding singularity-free saddle connection crosses its image. However, an important complication arises when one is careful about not conflating the rectangle with its projection to $S$; indeed, the argument really works in the universal cover and it applies to some lifted homeomorphism $\tilde{f}$. For this reason, it is conceivable that some of the fixed points one finds during the course of the proof actually project to the same point on the surface. So in the end, one obtains upper and lower bounds on the number of fixed points in terms of the number of times that the saddle connection crosses its image, \textit{and} the multiplicity of the projection map from the rectangle to its (potentially immersed) image on the surface. 

\item The singularity-free saddle connections can be organized into a so-called \textit{veering triangulation}\index{veering triangulation} of the mapping torus $M_{f}$ for $f^{\circ}$, the map obtained from $f$ by deleting all of its singularities. We refer the reader to Minsky--Taylor \cite{MT} and Gueritaud \cite{Gueritaud} for the details of this construction. For our purposes, it suffices to say that the veering triangulation is a topological ideal triangulation of the fibered manifold, and the edges of the $1$-skeleton of the  correspond in a precise way to singularity-free saddle connections. Moreover, as a set of arcs, these edges represent a totally geodesic subset in the arc and curve complex of the (punctured) surface. It will also suffice to think of the veering triangulation as belonging to the infinite cyclic cover $\widetilde{M_{f^{\circ}}}$. 

Minsky-Taylor \cite{MT} show that the triangulation encodes much of the same data as the model manifold for $\widetilde{M_{f^{\circ}}}$. Given a subsurface $Y \subset S$ for which $d_{Y}(\lambda^{-}, \lambda^{+})$ is sufficiently large, there is a sub-complex in the veering triangulation that is naturally homeomorphic to $Y \times [0,1]$. The $1$-skeleton of the base of this sub-complex is nearby (in the arc and curve complex of $Y$) to (the projection to $Y$ of) $\lambda^{-}$, and the arcs comprising the $1$-skeleton of the time $1$ slice are nearby to $\lambda^{+}$. Finally, there is a natural way to \textit{flip} through the triangulation by elementary moves. In this sense, one can interpret the veering triangulation as providing a combinatorial interpolation between the ending laminations for $f$: flipping through the triangulation in the upwards direction moves away from $\lambda^{-}$ and towards $\lambda^{+}$ in the sense that the higher up one moves, the more subsurfaces have the property that a projection of the triangulation on a given section resembles the projection of $\lambda^{+}$ to that subsurface. 

\item One can prove the existence of a section $T$ of the veering triangulation (which we call \textit{annular avoiding}) with the property that the projection map from the universal cover to $T$ has uniformly bounded multiplicity on each singularity-free rectangle. Thus, by $(1)$ above, the number of fixed points of $f$ can be related uniformly to the number of times that singularity free saddle connections cross their images in such a section. 

\item It remains to relate the intersection numbers between singularity-free saddle connections and their images, and the dilatation of the map $f$. It is here that we will need the key hypothesis that every simple closed curve intersects its $f$-image; using this (and some other ingredients from the theory of veering triangulations), given a section $T$ of the triangulation, one can relate the total intersection number $i(T, f(T))$, to the maximum value of $i(\sigma, f(\sigma))$ where $\sigma$ is a singularity-free saddle connection. 

In general, these two quantities can be arbitrarily far apart. Indeed, imagine a scenario in which every saddle connection $\sigma$ in the triangulation $T$ has the property that its image $f(\sigma)$ intersects some \textit{other} saddle connection $\sigma'$ many times, and maybe also imagine that $i(\sigma, f(\sigma))= 0$ for each such $\sigma$. In this situation, $i(T, f(T))$ will be very large but the maximum value of $i(\sigma, f(\sigma))$ is $0$. It's not obvious that this scenario is prevented by our assumption on $f$, but it is, and it is \textit{at least} believable that it is: think about for example the case where $f$ has a single singularity. Then each saddle connection in $T$ starts and ends at that singularity and therefore represents a simple closed curve; our assumption implies that the $f$-image of that curve must cross itself, whence the same must be true for each saddle connection in $T$.  

\item Finally, we only need to connect $i(T, f(T))$ to the dilatation of $f$. For this, we prove a version of a coarse equality due to Choi--Rafi \cite{CR} connecting the intersection number between a pair of markings to various sums of subsurface projections.

\end{enumerate}

\section{Open questions}

We leave the reader with some interesting open questions. 

\begin{question} \label{bounds} (sharp bounds) As discussed in Section \ref{geometry from combinatorics}, there exists several results relating the geometry of a mapping torus to the dynamics of the monodromy map on a combinatorial complex that are \textbf{effective}: the multiplicative and additive errors are bounded by known and computable functions of the topology of the fiber (see \cite{FS}, \cite{Vi}, \cite{APT}, \cite{ATW} for examples). How sharp are these functions? Can one construct a sequence of examples to demonstrate that they must grow with the topology of the fiber? 
\end{question}

A natural place to begin investigating Question \ref{bounds} would be in a fibered face of a fixed fibered $3$-manifold $M$, because this gives one access to many monodromies defined on many different surfaces whose corresponding mapping tori have identical geometric properties. Work of Minsky-Taylor \cite{MT} on veering triangulations and subsurface projection gives evidence that, in certain instances, one might be able to obtain uniform control over additive and multiplicative errors as one varies over a fibered face. 

\begin{question} \label{free geometry} (free groups) Can one use geometry in a way that is analogous to Section \ref{geometry from combinatorics} to study the dynamics of a fully irreducible element of the outer automorphism group of a free group? 
\end{question}

The difficulty inherent in Question \ref{free geometry} is that when it comes to a fully irreducible outer automorphism, there is no natural analog of the fibered hyperbolic $3$-manifold for which that automorphism serves as a monodromy. Nevertheless, Clay \cite{Clay} proved a very interesting analog of Kojima--McShane's theorem \cite{KJ} mentioned in Section \ref{WP}. The result of Kojima--McShane was originally stated as a bound on the volume of $M_{\phi}$ in terms of the \textit{entropy} of $\phi$, a quantity that measures the growth rate of word length upon iteration\footnote{The connection to Teichm{\"u}ller space is that the entropy is none other than the logarithm of the dilatation, a.k.a. the translation length on $\mathcal{T}(S)$ with respect to the Teichm{\"u}ller metric.}. In turn, Clay relates the entropy of an irreducible outer automorphism to its so-called $\ell^{2}$-torsion, which, by a result of L{\"u}ck and Schick \cite{LS}, captures the volume of a fibered hyperbolic manifold $M_{\phi}$. 

In search of results of this flavor, the author jointly with Clay and Rieck explored the Weil-Petersson like metrics mentioned in Section \ref{WP} that can be defined on the Culler-Vogtmann Outer Space\index{outer space} \cite{ACR}. The hope was that translation length with respect to this metric would also capture something like entropy or $\ell^{2}$-torsion. Unfortunately, we prove that the associated metric completion of the Outer Space has a global $\mbox{Out}(F)$-fixed point, and we amass evidence to suggest that the entire Outer Space may in fact have finite diameter with respect to these metrics. One can still search for a more combinatorial way to capture entropy, in the spirit of Theorem \ref{thm:brock1}:

\begin{question} (pants graph for the free group) Is there a ``pants graph'' for $\mbox{Out}(F_{n})$? In other words, is there some complex, supporting a simplicial action by $\mbox{Out}(F_{n})$, such that the translation length of a fully irreducible outer automorphism is coarsely equal to (the negative of\footnote{The convention is that this invariant is negative.}) $\ell^{2}$-torsion?
\end{question}

Aiming to generalize the results in Section \ref{geometry from combinatorics} in a different direction, one could consider whether there exist analogs in the setting of hierarchically hyperbolic groups:

\begin{question} (HHG analogs) Is there an ``HHG'' version of any of the theorems from Section \ref{geometry from combinatorics}? 
\end{question}

For example, the pants graph is itself a hierarchically hyperbolic space (HHS), built from factor spaces associated with curve complexes of all non-annular subsurfaces. One could therefore imagine interpreting it as an HHS associated to $\mbox{Mod}(S)$, in the sense that the latter acts simplicially on the former and the former is coarsely obtained from the latter by coning off all of the annular factors. Perhaps there is a more general statement lurking here, involving a hierarchically hyperbolic group and the HHS obtained from it by coning off the factors at the top of the hierarchy.

Heading back to the world of hyperbolic $3$-manifolds that fiber over the circle, we observe that the results summarized above revolve around volume, lengths of (shortest) closed geodesics, and electric circumference. There are of course many other interesting geometric features of $M_{\phi}$ that one might hope to relate directly to the action of $\phi$ on a combinatorial complex. For example, Baik--Gekhtman--Hamenst{\"a}dt estimate the smallest eigenvalue $\lambda_{1}(M)$ of the Laplacian in terms of the volume \cite{BGH}; it of course follows that there is some indirect relationship between $\lambda_{1}$ and $\tau_{\mathcal{P}}(\phi)$ but one could ask for a clearer connection between eigenvalues and combinatorics: 

\begin{question} Can one estimate higher eigenvalues of the Laplacian operator on $M_{\phi}$ in terms of combinatorial invariants associated with $\phi$?
\end{question}

A topic that we haven't explored at all in this piece is the so-called \textit{point-pushing} pseudo-Anosovs-- those pseudo-Anosovs that lie in the kernel of the homomorphism 
\[ \mbox{Mod}(S, p) \rightarrow \mbox{Mod}(S)\]
associated with forgetting the puncture, $p$. This kernel is naturally associated with $\pi_{1}(S,p)$ via the Birman exact sequence 
\[ 1 \rightarrow \pi_{1}(S,p) \xrightarrow{\iota} \mbox{Mod}(S,p) \rightarrow \mbox{Mod}(S) \rightarrow 1, \]
and a theorem of Kra \cite{Kr} implies that given $\gamma \in \pi_{1}(S,p)$, $\iota(\gamma)$ is pseudo-Anosov precisely when $\gamma$ fills $S$ in the sense that when one chooses a minimal position realization of $\gamma$, every complementary region to it is simply connected. 

There is a plethora of results that relate the combinatorics of a closed curve to various geometric properties of its geodesic realizations in some hyperbolic metric (see for instance \cite{AGPS}, \cite{Bas}, \cite{Bas2}, \cite{ArC}). It is therefore natural to wonder whether one can relate the dynamical properties of a point-pushing pseudo-Anosov map $\iota(\gamma)$ to the combinatorics of the associated closed curve, $\gamma$. 

To the best of the author's knowledge, very little is known on this front. Perhaps the most robust result to date in this area is due to Dowdall \cite{Dow} who relates the geometric self-intersection number $i(\gamma, \gamma)$ to the dilatation $\mbox{dil}(\iota(\gamma))$:
\[ (i(\gamma, \gamma) +1 )^{1/5} \leq \mbox{dil}(\iota(\gamma)) \leq 9^{i(\gamma, \gamma)}.\]

In joint work with Gaster \cite{AG}, we show that these bounds can be tightened for a \textit{random} point-push, i.e., a point-push associated to the destination of a very long random walk on the surface group. Indeed, letting $w_{n}$ denote this destination, we show that the probability that $\mbox{dil}(\iota(w_{n}))$ is bounded above and below by functions of the form $\exp(\sqrt{i(\gamma, \gamma)})$, converges to $1$ exponentially fast in $n$. However, tighter results that hold for a general point-push remain elusive. 

Beyond dilatation, there is a rich variety of other invariants one might hope to connect to the combinatorics of $\gamma$, for example any other geometric property of the mapping torus $M_{\iota(\gamma)}$, or the location of the axis of $\iota(\gamma)$ in the Teichm{\"u}ller space $T(S,p)$. We record this family of questions as follows: 

\begin{question} (point-pushes) What dynamical properties of a point-push pseudo-Anosov $\iota(\gamma)$ can be read off from the combinatorics of the filling curve $\gamma$?
\end{question}

\end{document}